\documentclass[reqno]{amsart}% 

\usepackage[utf8]{inputenc}
\usepackage[T1]{fontenc}
\usepackage[english]{babel}
\usepackage{amsmath}
\usepackage{amssymb}
\usepackage{amsthm}
\usepackage{systeme}
\usepackage{mathtools}

\numberwithin{equation}{section}
\usepackage[pagewise]{lineno}%\linenumbers 
\usepackage{amsaddr}%

\usepackage{enumerate} %opcional
\usepackage{graphicx}  %opcional
\usepackage{subcaption}
\usepackage[shortlabels]{enumitem}
\usepackage{bm} % 
\usepackage{emptypage} %
\usepackage{hyperref}
\usepackage{tikz-cd}

\usepackage{csquotes}

\usepackage[
    style=trad-alpha,
    maxalphanames=4
]{biblatex}
\usepackage{csquotes}

\usepackage{bm}
\newcommand{\vect}[1]{\boldsymbol{\mathbf{#1}}} %bold greek symbols and letters
\newcommand{\dd}{\textnormal{d}}

\usepackage{soul}
\usepackage[normalem]{ulem} %\sout allows to cross \cite. 

\newcommand{\sbullet}{\,\begin{picture}(-1,0)(-1,-1)\circle*{2}\end{picture}\ } %undfortunately radius can not change just a bit
\newcommand{\primebullet}[1]{\overset{\sbullet}{#1}}

\DeclareMathOperator{\tr}{tr}
\DeclareMathOperator{\sad}{sad}
\DeclareMathOperator{\cen}{cen}
\DeclareMathOperator{\midd}{mid}

\calclayout

\begin{document}
\title[Example simplest bifurcation diagram on a torus]{Example of simplest bifurcation diagram for a monotone family of vector fields on a torus}

\author{Claude Baesens, Marc Homs-Dones and Robert S. MacKay}
\address{Mathematics Institute, University of Warwick, UK.}
\thanks{M. H.-D. was supported by “la Caixa” Foundation (ID 100010434) with fellowship
code LCF/BQ/EU20/11810061, and by the Warwick Mathematics Institute Centre for Doctoral Training with funding from the University of Warwick. }
\email{marc.homs-dones@warwick.ac.uk}

\begin{abstract}
We present an example of a monotone two-parameter family of vector fields on a torus  whose bifurcation diagram we demonstrate to be in the class of ``simplest'' diagrams proposed by Baesens \& MacKay (2018 Nonlinearity 31 2928--81).  This shows that the proposed class is realisable.

%OLD: We show that the bifurcation diagram proposed in Baesens \emph{et al} (2018 \emph{Nonlinearity} \textbf{31} 2928-2981) is the simplest bifurcation diagram for a monotone two-parameter family of vector fields on a torus. We achieve this by finding an example that satisfies all the minimal conditions derived there. We give proofs for all but one of the conditions and present strong numerical evidence for the remaining one. 
%For the particular condition, existence of at most one contractible periodic orbit, we are only able to give strong numerical evidence. 
\end{abstract}



\date{\today}

\maketitle

\section{Introduction}

\subsection{Motivation}

It is now commonplace to compute bifurcation diagrams for families of dynamical systems, with one or more  parameters, e.g.~\cite{vanVeen03}. Very complicated diagrams can emerge, even for 2D ODE systems where no chaos is possible e.g.~\cite{Krauskopf94}.

The question arises whether there are topological criteria that force at least some amount of structure on the bifurcation diagram?  And if so, what are the simplest diagrams that can occur (where ``simplest'' will need defining)?

In the context of unimodal maps of the interval, for example, a one-parameter family is said to be {\em full} if it goes from a case with attracting fixed point to one for which the image of the critical point is a fixed point.  For a full family, an intricate sequence of bifurcations (including the famous period-doubling sequence) has to occur, governed by ``kneading theory''.  If the ``kneading sequence'' is monotone in the parameter then one obtains a simplest diagram, in the sense that for every other full family the kneading sequence takes some not necessarily monotone path so some bifurcations may be made and then undone and then remade.

In continuous time, an example was derived by Guckenheimer \cite{Guck1977}.  He showed that under a weak condition on a two-parameter family of gradient vector fields, the bifurcation diagram for equilibria has at least four cusps.  A simplest case (in the sense of minimising the number of cusps) is provided by Zeeman's catastrophe machine, for which the bifurcation set is precisely a closed curve with four cusps, called an astroid (many other examples giving an astroid were reviewed by Chillingworth \cite{Chil1987}).

\begin{figure}[htbp]
    \centering
        \includegraphics[width=1\textwidth]{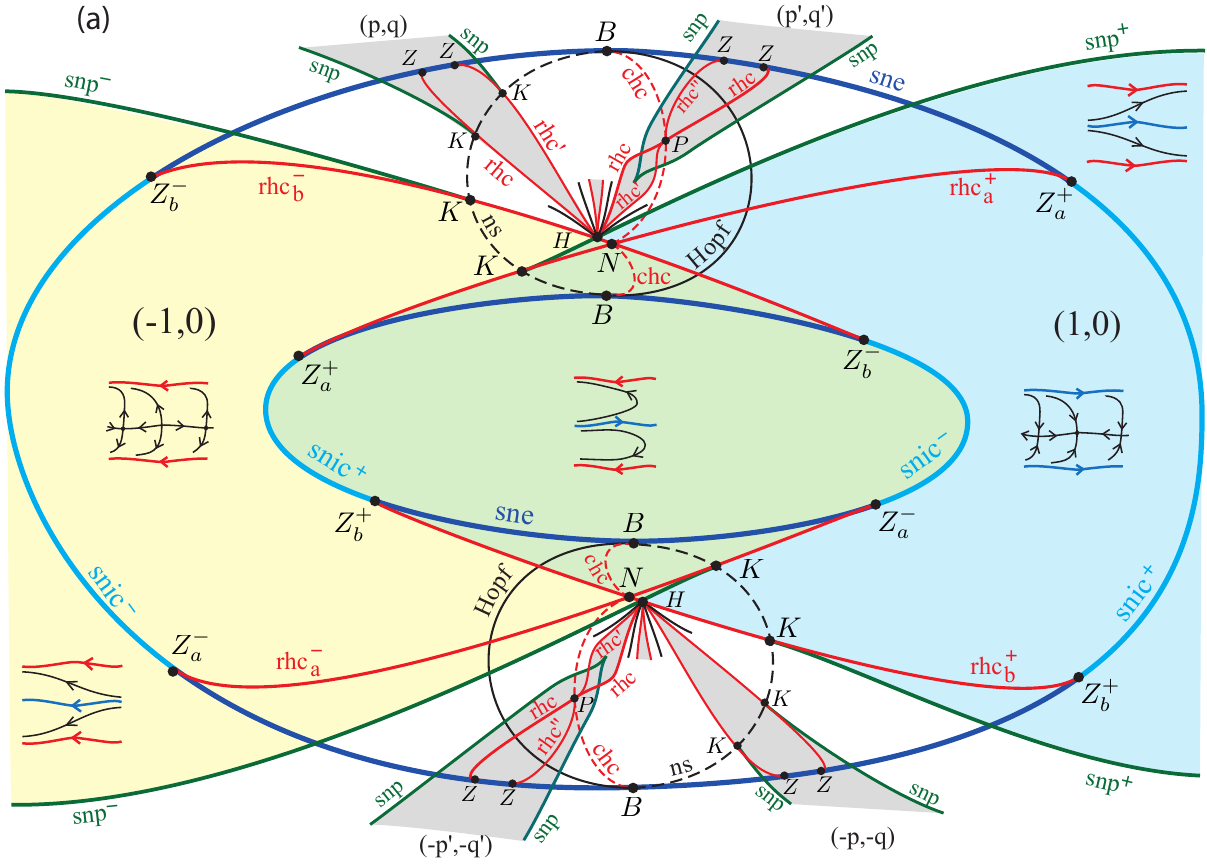}\\
        \vspace{8pt}
        \includegraphics[width=1\textwidth]{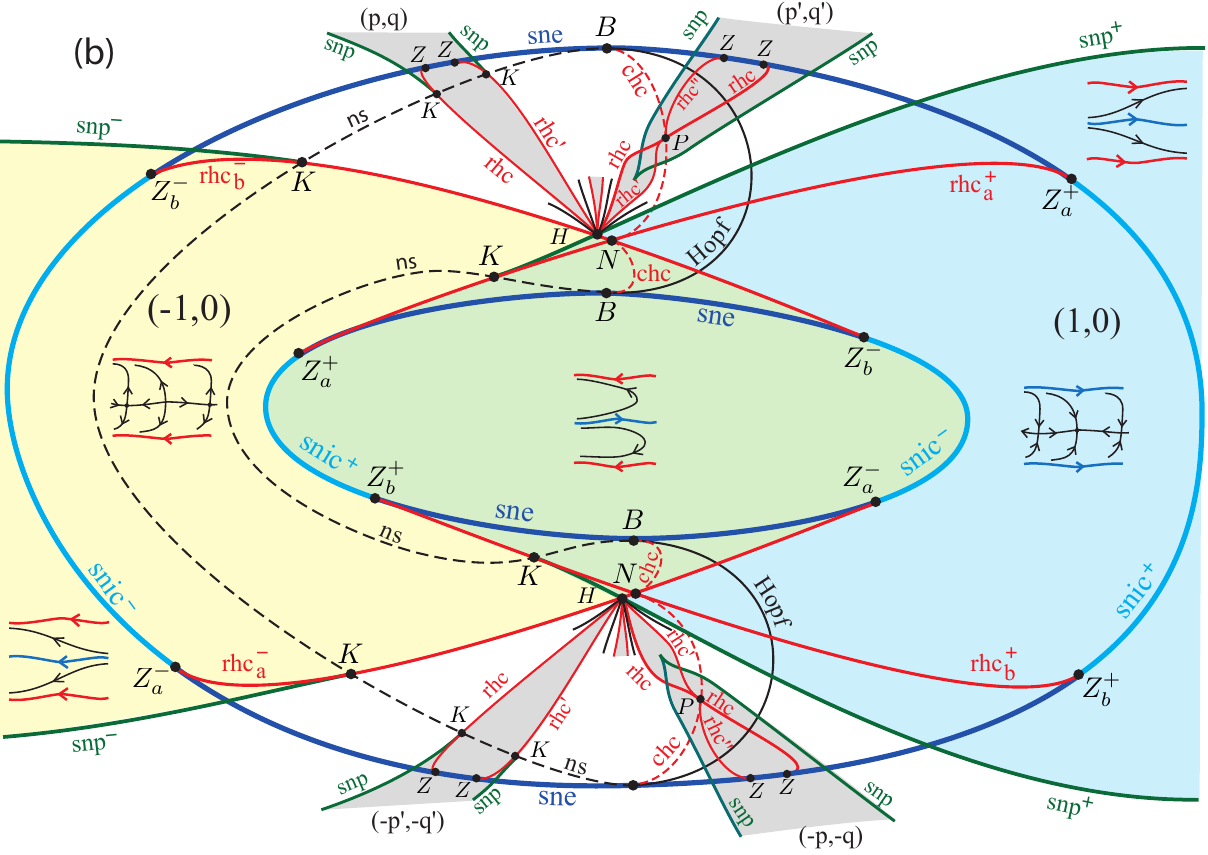}
    \caption{ Two simplest bifurcation diagrams for monotone families of vector fields on a torus  (slightly edited version of \cite[Figure 22]{Baesens_2018}). In this work we will focus on (a). The infinitely many rhc curves and %corresponding 
    tongues emanating from the $H$ point have not been fully depicted.   }
    \label{fig:full_bif_diag}
\end{figure}

More recently, this philosophy was applied to 2-parameter families of vector fields on a 2-torus.  Under a ``monotonicity'' assumption on the family, it was shown in \cite{Baesens_2018} that the bifurcation diagram must be at least as complicated as one of those in figure~\ref{fig:full_bif_diag}, in a sense that we shall recall.  

An example was given there (a mistake in its analysis was corrected in \cite{Baesens_2022_Corrigendum}) for which the bifurcation diagram is only slightly more complicated than 
figure~\ref{fig:full_bif_diag}(a):~it contains four points of degenerate Hopf bifurcation and two regions with two contractible periodic orbits, whereas figure~\ref{fig:full_bif_diag}(a) has no degenerate Hopf points and at most one contractible periodic orbit for given parameter values. So the question remained whether figure~\ref{fig:full_bif_diag}(a) is realisable.
 
In the present paper, we study a modification of the example, prove that it still satisfies the conditions that were proved for the example in \cite{Baesens_2018}, that it has no degenerate Hopf points, and show numerically that it has at most one contractible periodic orbit.  Thus the case of figure~\ref{fig:full_bif_diag}(a) is indeed realisable and therefore simplest. 

\subsection{Background and description of the bifurcation diagram}

We now proceed to give a detailed explanation of the simplest bifurcation diagram depicted in figure~\ref{fig:full_bif_diag}(a). The rest of this work will focus on showing that the forthcoming family \eqref{eq:principal} has this type of bifurcation diagram. 

The main division in figure~\ref{fig:full_bif_diag}(a) is given by the  curves of parameter values where the vector field has \emph{saddle-node equilibria} (sne), 
%EXPAND: not defined regorously purposly as we are a bit fusi on weather B points are included in this curves or not (depending on the context we change). If we have to be rigorous probably best they are included but not in the generic saddle nodes.
which are represented by two tones of blue and form two topological circles. The closed annulus formed by the parameter values between the circles is the \emph{resonance region}, which corresponds to parameters where the vector field has equilibria. The complement of the resonance region is divided into two components, one bounded which we refer to as the \emph{hole}, and an unbounded one named the \emph{outside}. 
The saddle-node curves have arcs in light blue which correspond to parameters with a \emph{saddle-node on an invariant circle} (snic). These invariant circles are in fact \emph{rotational} meaning homotopically non-trivial. The red solid curves correspond to parameter values with a \emph{rotational homoclinic connection} (rhc),  which is a branch of stable manifold of a saddle connecting to a branch of unstable manifold of the same saddle but making a rotational loop.
We can see in figure~\ref{fig:full_bif_diag}(a) that the rhc curves meet the sne ones in a $Z$ \emph{point}, a parameter value where a saddle-node equilibrium has a rotational homoclinic connection between its neutral direction and its hyperbolic one. A \emph{necklace point}
($N$) is a parameter point at which there are rhcs in opposite directions connecting the same saddle. Emanating from it we can see two parameter curves of \emph{contractible homoclinic connection} (chc) in dashed red. 
% llll DEFINE P point? I don't think there is a need as not discussed in this work.

The other main division in figure~\ref{fig:full_bif_diag}(a) is given by the shaded areas. The blue and yellow areas correspond to the existence of periodic orbits or invariant cycles of \emph{horizontal} homotopy type, meaning $(\pm 1,0)$, as we can see depicted in the phase diagrams. We note that in the hole the phase diagram consists of two \emph{Reeb components}, i.e.~annuli bounded by periodic orbits in opposite directions and containing no equilibria nor other periodic orbits. 
Parts of the boundary of the shaded regions consist of curves of parameter values with  rhc, where a rotational periodic orbit is generated or destroyed as the connection is broken. The rest are curves of \emph{rotational saddle-node periodic orbit} (snp), i.e.~a periodic orbit that is attracting from one side and repelling from the other, where two periodic orbits with opposite stability coalesce and annihilate each other. 
% EXPAND: that is what happens generically 
%old: We see that the boundary of the shaded regions are  parameter curves of either rhc or rotational \emph{saddle-node periodic orbit} (snp), that is a periodic orbit that is attracting from one side and repelling from the other. \textcolor{brown}{In the rhc a periodic orbit becomes the connection which is then destroyed, whereas in a snp two periodic orbits merge into a snp one which then disappear}.
In the upper and lower part of the diagram we can see grey shaded tongues which correspond to other homotopy types.
% EXPAND: They correspond to rational winding reatio, but we don't get into this in this work. Note that homotopy type by definition must be rational winding ratio.  
Although not represented fully, there are in fact infinitely many of them emanating from an $H$ \emph{point}, which is a parameter value where we have coexistence of a rhc with a snp in opposite directions. The tongues that have been fully depicted together with their corresponding bifurcations (rhc, snp curves and $Z$, $K$, $P$ points) only give an example of how they may extend from the $H$ point and are not necessary conditions for simplicity. 

In solid black we can see a parameter curve where a generic Hopf bifurcation occurs, so for these parameter values the vector field has a \emph{centre}, i.e. an equilibrium with purely imaginary eigenvalues. In dashed black we represent curves where the vector field has a \emph{neutral saddle} (ns), meaning that its eigenvalues are equal with opposite sign. The Hopf and ns curves meet in $B$ \emph{points}, parameter values having an equilibrium with double eigenvalue 0 and non-diagonalisable, and in our example a generic Bogdanov–Takens bifurcation occurs. The curves of ns intersect the curves of rhc at  $K$ \emph{points}, parameter values where the flow has a neutral saddle with a rhc. To complete our terminology, although it does not appear in the figure, a $J$ \emph{point} is a parameter value where the vector field has a ns with a chc. 

The unfamiliar reader can find depictions of the generic local unfoldings of all the bifurcation points discussed above in \cite{Baesens_2018} (except for $J$ points but their unfolding is the same as for $K$ points with rotational replaced by contractible).
 % EXPAND:  There you can also find definition of P points. 

\subsection{Results}
The setting is two-parameter families of vector fields
\begin{equation}
    \dot {\vect x}= \vect G(\vect x, \vect \Omega)
    \label{eq:x'=G}
\end{equation}
with $\vect x\in \mathbb T^2= \mathbb R^2/\mathbb Z^2$, $\vect \Omega \in \mathbb R^2$ and $\vect G\in \mathcal C^3$. We assume them to have monotone dependence on $\vect \Omega$ in the sense that for a given inner product $\langle \cdot ,\cdot\rangle $ and $c>0$,
\[\langle D_{\vect \Omega} \vect  G _{({\vect x} ,{\vect \Omega} )}\: \vect w, \vect w \rangle \geq c|\vect w|^2\]
for all  $({\vect x} ,{\vect \Omega} )\in \mathbb T^2 \times \mathbb R^2$ and  $\vect w\in \mathbb R^2$, where $D_{\vect \Omega} \vect G _{({\vect x} ,{\vect \Omega} )}$ is the derivative with respect to the second argument at $({\vect x} , {\vect \Omega} )$.  We want our family to be generic, meaning that  the main bifurcations that occur in the family will correspond to codimension 1 or 2 degeneracies, and our family will be transversal to the corresponding submanifolds of degenerate vector fields. 
Thus our family belongs to a residual set, i.e.~a countable union of open dense sets in the space of $\mathcal C ^3$ families.

The monotone condition may seem quite arbitrary at first, but it is a natural generalization of the straightforward class of families 
\begin{equation}
\begin{split}
    \dot x &= \Omega_x + v(x,y)\\
    \dot y &= \Omega_y + w(x,y).
\end{split}
\label{eq:v_w}
\end{equation}
Moreover, these families are thematically appropriate as the parameter affects them by a simple addition. The generic condition also plays an important role here, as for example it rules out the trivial choice of $v, w = 0$, which yields the family of constant vector fields on the torus.   This family is well understood but we do not consider it simplest, as every parameter value has a bifurcation from  infinitely many periodic orbits to none, which is a highly degenerate phenomenon.
% Note: we need to  restric to a family, if we consider simples in any family, constant family (unchanged by changing parameters) is obviouly the simplest and it is uninteresting. 

We take the same approach as in \cite{Baesens_2018} to define ``simplest''. In that work it was first shown that any family satisfying the conditions below \eqref{eq:x'=G} has parameter values with at least two equilibria. So the first criterion for simplicity is to minimize the maximal number of equilibria, and the first assumption for the simplest bifurcation diagram is to have at most two. Continuing the study of the equilibria under this assumption it was found that the set of parameter values with equilibria is an annulus  bounded by  parameter curves of saddle-node, which were shown to contain at least four Bogdanov–Takens points. So the second simplicity criterion is to minimize the number of $B$ points, and the simplest will have at most four. Sequentially continuing in this manner, that is, finding a feature that complicates the bifurcation diagram and then minimizing the instances of it, a nine point definition of simplicity was devised in \cite{Baesens_2018}. Moreover, it was shown that a bifurcation diagram satisfying the following assumptions would be simplest.

\begin{enumerate}[leftmargin=19pt,labelsep=3pt]
%lll change left marigin if change format. 

\item[\textbf{1.}] there are at most two equilibria for any parameter values; 

\item[\textbf{2.}] there are at most four $B$ points; 

\item[\textbf{3a.}] there is no closed curve of centres nor closed curve of neutral saddles;

\item[\textbf{3b.}] there is no coexistence of a centre and a neutral saddle; 

\item[\textbf{4a.}] there are at most two Reeb components for parameter values in the hole;

\item[\textbf{4b.}] there is no other invariant annulus for flows in the hole;

\item[\textbf{4c.}] there is no Reeb component for parameters outside the resonance region; 

\item[\textbf{5a.}] there are at most four $Z$ points on the inner saddle-node curve;

\item[\textbf{5b.}] there are at most four $Z$ points of horizontal homotopy type on the outer saddle-node curve (it may have $Z$ points of other homotopy types);

\item[\textbf{5c.}] there are at most four curves of rotational homoclinic connection (rhc) of horizontal homotopy type; 

\item[\textbf{6.}] there are at most two necklace ($N$) points; 

\item[\textbf{7a.}] there is no $J$ point;

\item[\textbf{7b.}] there is at most one contractible periodic orbit (cpo) for any parameter value; 

\item[\textbf{8.}] there are at most four $K$ points of horizontal homotopy type; 

\item[\textbf{9.}] there are at most two $H$ points.

\end{enumerate}

We note that sequentially minimizing instances of features will give rise to a total preorder between bifurcation diagrams, and then by simplest we really mean a minimal element. So it is in principle possible to have more than one simplest bifurcation diagram. However, once we show that the assumptions above are realisable by attaining figure~\ref{fig:full_bif_diag}(a), then by  \cite{Baesens_2018} the only other possible simplest bifurcation diagram is given by figure~\ref{fig:full_bif_diag}(b), which would also satisfy them.  

%old: (before bakground)\mh{ We note here that the relation of being ``simpler'' that we will define give rise to a preorder, and by simplest we really a mean minimal element. In particular, it is possible to have more than one simplest bifurcation diagram, as is shown by the possibility of figure~\ref{fig:full_bif_diag}(b).}

We prove (sections \ref{sec:first_prop}, \ref{sec:delta_neigh} and \ref{sec:K_H_points})  that the family of vector fields, 
\begin{equation}
\begin{split}
    \dot x &= \Omega_x-\cos2\pi (y-\phi)-\varepsilon \cos2\pi x\\
    \dot y &= \Omega_y-\sin 2\pi y-\varepsilon \sin 2\pi x
\end{split}
\label{eq:principal}
\tag{\#}
\end{equation}
where $\phi\in (\frac{1}{24},\frac{5}{24})$ and $\varepsilon$ is small and positive,
satisfies all the above assumptions except 7b, and present strong numerical evidence that it satisfies 7b (section \ref{sec:cpo}). To obtain a formal proof one would need an insight or rigorous computer-assisted estimates (in the vein of \cite{computer-assisted}) to determine the sign of the integrals \eqref{eq:tilde_alpha_prime}, which would automatically give at most one cpo. Note that the proposed family is in the class of \eqref{eq:v_w}, so that it also realises the simplest bifurcation diagram in that class.
For future reference note that  
\begin{equation}
    \cos 2\pi (y-\phi)=C\cos2\pi y + S \sin 2\pi y,
    \label{eq:C_S_expression}
\end{equation}
% and,
% \begin{equation*}
%     \sin 2\pi (y-\phi)=C\sin 2\pi y - S \cos 2\pi y,
% \end{equation*}
where $C=\cos 2\pi\phi>0$, $S=\sin 2\pi\phi>0$.

In \cite{Baesens_2018}, the case $\phi=0$ was considered, but it was found that there were four points of degenerate Hopf bifurcation, which caused the appearance of more than one cpo in some regions of parameter space, thus contradicting assumption 7b. This motivated our introduction of $\phi$, which enables us to tune the system, and as we will see in appendix \ref{append:lyapunov_coef}, eliminates the degenerate Hopf points when $\phi\in (\frac{1}{24},\frac{5}{24})$. Another apparent obstruction to satisfying assumption 7b was that in \cite{Baesens_2018} it was erroneously computed that the necklace points were outside the trace-zero loops. This error was corrected in the corrigendum \cite{Baesens_2022_Corrigendum} together with numerical computation of the relevant part of the bifurcation diagram.

Unavoidably in this work we have repeated many of the figures and arguments from \cite{Baesens_2018}. Some of them have not changed with the introduction of $\phi$ whereas others have become significantly more intricate. We have also taken the opportunity to flesh out arguments that were only sketched there, for instance see the study of snp curves in section \ref{sec:transit_map} and of chc curves in section \ref{sec:lack_of_chc_cpo}.

\section{First properties }
\label{sec:first_prop}

Besides the assumptions listed in the introduction, we also need to check that our family is $\mathcal C^3$ and monotone, both of which are immediate for our example (considering the standard inner product). Moreover, for the principal features of the bifurcation diagram (saddle-node equilibria, Hopf bifurcation, chc, horizontal rhc and $B$, $Z$, $N$, $K$, $H$ points),
we will check throughout the text that they are codimension at most two with our family being transverse to them. It is probably unfeasible to check this for all features, as in particular it would require a study of the infinitely many curves of rhc emanating from the $H$ points, but  appealing to $\mathcal C^3$-genericity we can deduce that there are arbitrary small $\mathcal C^3$ perturbations of our example satisfying this.

\subsection{Equilibria}
\begin{figure}[htbp]
     \centering
     \begin{subfigure}[b]{0.47\textwidth}
         \centering
         \includegraphics[width=\textwidth]{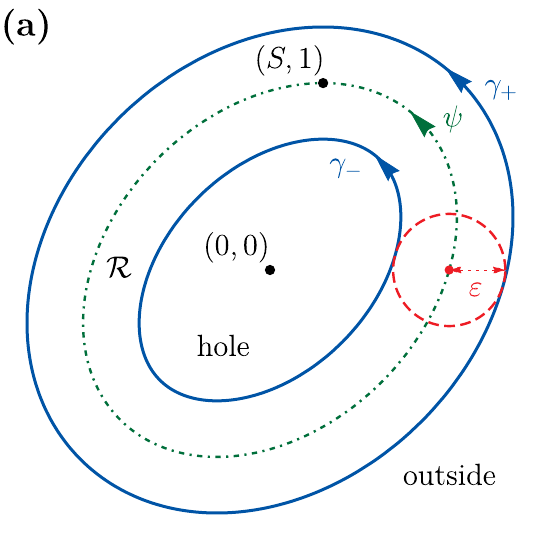}
     \end{subfigure}
     \hfill
     \begin{subfigure}[b]{0.47\textwidth}
         \centering
         \includegraphics[width=\textwidth]{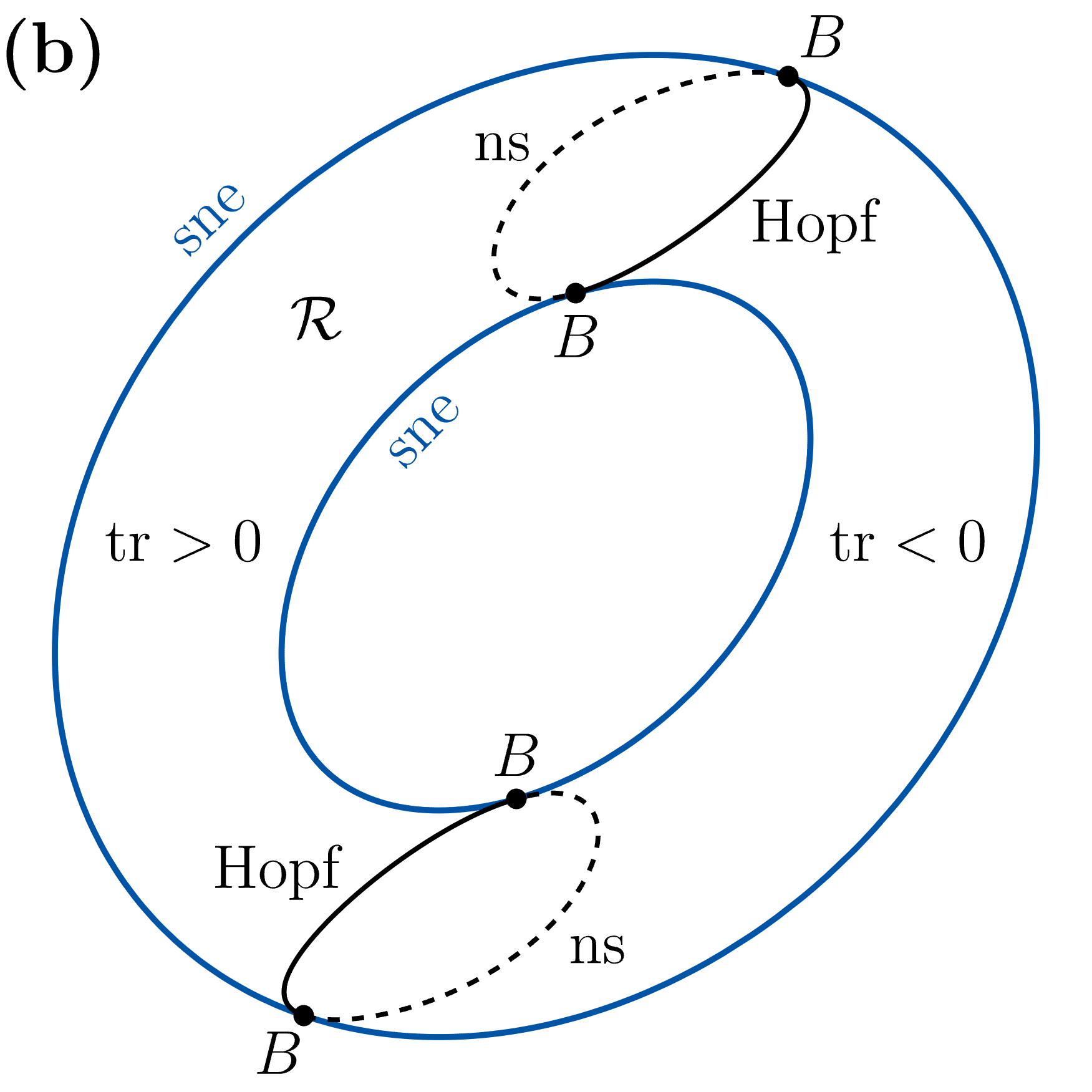}
     \end{subfigure}
        \caption{(a) The resonance region $\mathcal R$; (b) Bifurcation diagram for the equilibria. To depict these figures, as well as the upcoming ones, we have chosen $\phi=1.1/24$ for aesthetic reasons.  }
        \label{fig:resonance}
\end{figure}

The equilibria of \eqref{eq:principal} are given by,
\begin{equation}
\begin{split}
     \Omega_x &=\cos2\pi (y-\phi)+\varepsilon \cos2\pi x\\
    \Omega_y &= \sin 2\pi y+\varepsilon \sin 2\pi x.
\end{split}
\label{eq:equilibria}
\end{equation}
The set $\mathcal E$ of equilibria $(\vect \Omega, \vect  x)$ as a subset of $\mathbb R ^2\times \mathbb T^2$ is the graph of the function $\vect \Omega ^e: \mathbb T^2 \rightarrow \mathbb R^2$ defined by \eqref{eq:equilibria}, and thus it is a topological torus.
%\expand{It is a general fact. Think how a function in the circle would be just a deformed circle, same is true for torus (or any surface right?).}
The resonance region $\mathcal R$ is the projection of $\mathcal E$ onto $\mathbb R^2$ or equivalently the image of $\vect \Omega ^e$.
%\expand{As R is parameters where we have equilibria. Thus is first component of E.}
Furthermore, when $y$ and $x$ rotate around the torus, $\vect \Omega^e(\vect x)$ performs similar movement to an epicycle, where the base circle is replaced by the ellipse\footnote{Note that when $\phi =\pm \frac{1}{4}$, the ellipse $\psi$ degenerates to a line segment. This does not happen in our context as $\phi \in (\frac{1}{24},\frac{5}{24})$.} \begin{equation}
    \psi(y)=(\cos 2 \pi (y-\phi), \sin 2\pi y)
    \label{eq:psi}
\end{equation}
and the other is a circle of radius $\varepsilon$, which is depicted in dashed red in figure~\ref{fig:resonance}(a). 
%\expand{y determines position in elipse. x determines position in the circle centered at the point in the elipse. }
In particular, the resonance region is a closed $\varepsilon$-neighbourhood of the ellipse $ \psi$.
%\expand{It is clearly contained in a $\varepsilon$-neighobourhood, as all points are at most $\varepsilon$ far away from it (being in a cercle of this radius centered in a point in elipse). When moving y and thinking of what the circle around the point in the elipse covers, it is clear it toches all points in neighbourhood. }

%IMPROVEMENTS:  move phi and gamma in the first figure (so not so close to arrow). Maybe also improve this arros that can be abit unclear). 
% Maybe also increase slightly font and size points although things may start getting tight (and font now is enoguh to be comparable with text, so should be ok)

\subsubsection{Parametrisation of the resonance region}
\label{sec:param_resonance}
For $\varepsilon$ small enough we can parametrise the $\varepsilon$-neighbourhood of $\psi$ by its tubular neighbourhood, i.e. $\mathcal R$ can be parametrised by,
\begin{equation}
    \vect \Omega = \gamma_r(\theta):= \psi (\theta) + r N_\psi (\theta), 
    \label{eq:ellipse_parametrisation}
\end{equation}
where $N_\psi$ is the external unit normal vector of $\psi$, $\theta \in \mathbb T^1$ and $|r|\leq\varepsilon$. So in this case, the boundary of $\mathcal R$ is given by 
\begin{equation}\gamma_\pm (\theta):=\gamma_{\pm \varepsilon}(\theta)= \psi (\theta)\pm \varepsilon N_\psi (\theta).
 \label{eq:def_gamma}
\end{equation}
By studying for which $r$  the curve $\gamma_r$ stops being regular one finds that \eqref{eq:ellipse_parametrisation} is injective, and thus
gives a parametrisation,  as long as $r>-R_{\min}$ where $R_{\min}$ is the minimal radius of curvature of $\psi$.
So $\varepsilon$ small enough in the start of this section can be quantified by $\varepsilon<R_{\min}$. Moreover, for $r>\varepsilon$, equation \eqref{eq:ellipse_parametrisation} is also a parametrisation that covers the outside of $\mathcal R$.

\subsubsection{Checking assumption 1}
We check that when  $\varepsilon <R_{\min}$, there are precisely two equilibria for each interior point of $\mathcal R$ and one for each boundary point of $\mathcal R$, so that assumption 1 is satisfied. 

For a fixed $y$ and varying $x\in \mathbb T^1$, $\Omega^e(\vect x)$ makes a circle centred at $\psi(y)$ of radius $\varepsilon$, which by \eqref{eq:def_gamma} intersects the inner (resp.~outer) boundary of $\mathcal R$ at precisely one point, see red dashed circle in figure~\ref{fig:resonance}(a).  Then, doing a positive lap of  $y$  moves this circle and the intersections anticlockwise %, as $\phi \in (-\frac{1}{4},\frac{1}{4})$,
 making a full lap around the hole of  $\mathcal R$.
Thus, each of the two arcs of the circle with end points in $\partial \mathcal R$ goes through every point $\vect \Omega \in \mathcal R$ precisely once. Hence, we have exactly two equilibria for parameters in the interior of $\mathcal R$ and only one in the boundary as both arcs coincide there. 
%EXPAND: An alternative direct proof is given by taking point where we have atleast 3 roots then making circle radius epsilon. Then can show at some point boundary we have atleast 2 solutions, which implies no tubular neighbourhood  (as in boundary (cos 2pix,sin 2pix) vector is  normal to psi) so contradiction. 

\subsection{Linearisation at equilibria}

Given parameters $\vect \Omega$ let us linearise the vectorfield around an equilibrium $(x,y)$, for tangent orbit $(\delta x, \delta y)$, 
\begin{equation}
    \begin{pmatrix}\dot {\delta x} \\\dot {\delta y}\end{pmatrix}=
    \begin{pmatrix}
    2\pi \varepsilon \sin 2\pi x & 2\pi \sin 2\pi (y-\phi) \\
    -2\pi \varepsilon \cos 2\pi x & -2\pi \cos 2\pi y
    \end{pmatrix}
      \begin{pmatrix}{\delta x} \\{\delta y}\end{pmatrix}.
      \label{eq:linearization_equilibria}
\end{equation}
The matrix determinant can be written as
\[\det =- 2\pi \varepsilon \left \langle(\cos 2\pi x, \sin 2\pi x) \hspace{1pt} , \hspace{1pt} \psi'(y) \right \rangle,\]
where $\psi$ was defined in \eqref{eq:psi}, so $\det =0$ when $(\cos 2\pi x, \sin 2\pi x)\bot \psi'(y)$, or equivalently when $(\cos 2\pi x, \sin 2\pi x)=\pm N_\psi (y)$. By \eqref{eq:def_gamma}, this happens precisely when $\vect \Omega$ is in the boundary of $\mathcal R$. 
%\expand{Simply note that here we are in det of equilibria, and the parameters corresponding to this equilibria existing correspond to $y=\theta$ so that using \eqref{eq:def_gamma} makes sense. }
So in the boundary of $\mathcal R$ we have saddle-nodes.  As shown in appendix \ref{append:generic_saddle-node}, they are generic saddle-nodes except when the trace vanishes, which instead gives generic $B$ points as discussed in the next subsection. 
Moreover, the torus $\mathcal E$ can be viewed as two copies of $\mathcal R$  with determinant of opposite signs,  glued together along the $\det = 0$ curves. So for each parameter value in the interior of $\mathcal R$, one equilibrium has Poincaré index +1 and the other -1, which is a saddle.

The linearised vector field has trace, 
\[\tr = 2\pi (\varepsilon\sin 2\pi x- \cos 2\pi y).\]
Thus the trace is zero where $\cos 2\pi y=\varepsilon\sin 2\pi x $, which forms two non-contractible closed curves on $\mathcal  E$, one near $y=\frac{1}{4}$ and the other one near $y=-\frac{1}{4}$, dividing $\mathcal  E$ into two pieces.  The trace of both equilibria is negative on the part projecting to the right of these curves, and positive to the left, whereas in the region enclosed by the projection of these curves (we will show they are simple curves in section \ref{sec:trace-zero_curves}) the sign differs between equilibria, see  figure~\ref{fig:resonance}(b).

\subsubsection{\texorpdfstring{$B$}{B} points}
We study when  both $\tr$ and $\det$ vanish. Let $x(y)$
%EXPAND: note that as x(y) is rotation normal elipse actual x(y) is monotone (increasing) but I think not needed as (x(y), y) is a a simple curve anyways due to y
be implicitly defined by one of the curves $\det=0$ and recall that $\tr$ can only vanish at $y=\pm\frac{1}{4}+O(\varepsilon)$. Then
\[\frac{\dd}{\dd y}\tr (x(y),y)|_{y=\pm \frac{1}{4}+O(\varepsilon)}  =\pm 4\pi ^2 + O(\varepsilon).\]
%\expand{$tr=0$ iff $\varepsilon\sin 2\pi x= \cos 2\pi y$ so $y\sim \pm 1/4$. Now, $\cos 2\pi (\pm 1/4+a) = \mp \sin 2\pi a =\mp 2\pi a+ O(a^3)$, so $y=\pm \frac{1}{4}+O(\varepsilon)$.  Now to compute derivative we have,             \[\frac{\dd}{\dd y}\tr (x(y),y) =\frac{\partial}{\partial x}\tr \cdot x'(y) + \frac{\partial}{\partial y} tr = 4\pi^2(\sin 2\pi y + \varepsilon x'(y)\cos 2\pi x) \]      Evaluating we get the result. (Note: equation defining $det =0$ does not involve $\varepsilon$, so $x'(y)$ neither does and is $O(1)$). To conlcude that we have simple zero we simply note that in $\det =0$ the value of the trace is monotone in this region (checked by derivative) so have simple zero. The fact that we have atleast a 0 is clear as the trace in the right sided and left side of the det=0 curves differs.   }
So we have a single simple zero close to $y=\frac{1}{4}$ and another one close to $y=-\frac{1}{4}$ for each $\det=0$ curve. These 4 points have  double eigenvalue zero, with non-diagonal Jordan form as the entries of \eqref{eq:linearization_equilibria} can not all simultaneously vanish. Thus, they are $B$ points and  assumption 2 is satisfied. We check that they are generic $B$ points in appendix \ref{append:generic_B}.

\subsubsection{Trace-zero curves}
\label{sec:trace-zero_curves}
The part of each trace-zero curve with index +1 consists of centres, and the trace has non-zero derivative across it,
%\expand{see appendix (and exapnd in it). }
% EXPAND old: Changing slightly parameters corresponds to changing x,y slightly. Now grad tr=2pi (-epsilon sin 2pi x, sin 2pi y) and in tr=0 sin2pi y close to \pm 1, so non degenerate 0, so across curve not zero derivative (note is across level curve tr=0, so would need two directions with no increment to have across it 0 derivative, would mean gradient (0,0)). Note, we can still have intersections, as we project in R which can create this. Intersections correspon to to different equilibria with trace 0, but both of them would have generic bifurcations.  This gives transversality
so these are curves of Hopf bifurcations. Crucially, we show in  appendix \ref{append:lyapunov_coef}, that they are generic Hopf bifurcations, in contrast to the system with $\phi=0$, studied in \cite{Baesens_2018}, which had four points of degenerate Hopf bifurcation. Thus, it is possible that our example satisfies assumption 7b, as all its Hopf points generate only one cpo. 

The part of each trace-zero curve with index -1 consists of neutral saddle (ns) and the derivative of the trace across it is non-zero too, which will be needed later. 
%\expand{This is checked in appendix when studying Hopf bifurcation, as we check transversality for the whole tr=0 curve. }
So each curve of trace-zero consists of an arc of centre and an arc of ns separated by $B$ points as depicted in figure~\ref{fig:resonance}(b).
%\expand{This actually proofs assumption 3a if we interpret closed curve as "isolated" closed curve, as it was done in the origial paper. There, self intersection of ns or center was discarted as we only have 1 center or one ns so self intersections would be codimension 3. Here we allow for interpretation of closed curve without isolated (eg. and arch that has self intersection generates closed curve) so we leave this conclusion for after showing trace zero curve is simple. Note that this possible reinterpretation doesn't not change the assumptions, as when working in codimension at most 2 and with at most one saddle/center, as we are doing here, this conditions are equivalent. }
They have $y=\pm\frac{1}{4}+O(\varepsilon)$, so they are within $O(\varepsilon)$ of $\Omega_x=\pm S$, $\Omega_y=\pm 1$ but we can be more precise.
%\expand{We justify O term in parameters in the following equation (although can be seen directly more easely by ingnoring x terms (have $\varepsilon$) and doing taylor on y).}
 Using  that $\varepsilon\sin 2\pi x = \cos 2\pi y $ together with \eqref{eq:equilibria} and \eqref{eq:C_S_expression}, we find that the trace-zero curves are parametrised to first order by
 \begin{equation*}
\begin{split}
     \Omega_x &=\pm S+ \varepsilon(\cos2\pi x+C\sin 2\pi x)+ O(\varepsilon^2)\\
    \Omega_y &=\pm 1+\varepsilon \sin 2\pi x+O(\varepsilon^2),
\end{split}
\end{equation*}
%\expand{Substituding cosinus for \eqref{eq:C_S_expression} and doing taylor at $y=\pm\frac{1}{4}+O(\varepsilon)$ for sinus we get       \[\begin{split}\Omega_x &=C\cos2\pi y + S \sin 2\pi y+\varepsilon \cos2\pi x = C\varepsilon \sin 2\pi x+ S (\pm 1 +O(\varepsilon^2))+ \varepsilon\cos 2\pi x     \\     \Omega_y &= \sin 2\pi y+\varepsilon \sin 2\pi x = (\pm 1+O(\varepsilon^2))+\varepsilon \sin 2\pi x\end{split}\]         and the desired expression follows.  }
where the derivative with respect to $x$ of the error terms is also $O(\varepsilon ^2)$. 
So we deduce that the trace-zero curves are mapped to simple closed curves in the parameter space $\mathcal C^1$-close to the ellipses,
\begin{equation}
    \left (\Omega_x- C\Omega_y \pm (C-S)\right )^2+\left ( \Omega_y\mp 1\right )^2 = \varepsilon^2.
    \label{eq:trace_zero_elipse}
\end{equation}
%\expand{\[\begin{split} \Omega_x \mp S&= \varepsilon(\cos2\pi x+C\sin 2\pi x)+ O(\varepsilon^2)\\  \Omega_y \mp 1 &=\varepsilon \sin 2\pi x+O(\varepsilon^2),\end{split}\]      so      \[\begin{split} \Omega_x \mp S  -C(\Omega_y \mp 1)&= \varepsilon \cos2\pi x+ O(\varepsilon^2)\\  \Omega_y \mp 1 &=\varepsilon \sin 2\pi x+O(\varepsilon^2),\end{split}\]       from which the ellipse equation follows. \\ This implies that our curve is simple as derivative is also close to ellipse which is never close to 0 so can't have small intersection (then at some point we would have derivative very different in direction and thus very different in general as it is not close to 0). We can have far intersection as ellipse is simple.   } 
In particular there is no simultaneous centre and
ns, proving assumption 3b. As there are no  closed curves of centre nor closed curves of ns, assumption 3a is also satisfied. 
%\expand{As mentioned in an expand above (just after mentioning figure 2b), this could have been shown sooner. }

\subsection{Flow outside the resonance region}

Outside of $\mathcal R$, meaning in the unbounded component of the complement of $\mathcal R$,
we show there is Poincar\'e flow,  i.e. the flow has a global cross-section.
%EXPAND: (a cross-section for a flow is a transverse section such that the forward and backward trajectories of every point hit it.

Let $\vect \Omega $ be outside the resonance region, then by the comments made in section \ref{sec:param_resonance}, we have $\vect \Omega =\psi (\theta)+rN_\psi (\theta) $, for some angle $\theta$ and $r>\varepsilon$. Then, defining
\[h=\langle N_\psi(\theta), (x,y)\rangle\]
(on the universal cover $\mathbb R^2$ of $\mathbb T^2$), we have, 
\begin{equation*}
\dot h=\left \langle - N_\psi(\theta) \hspace{1pt},\hspace{1pt} \psi (y) - \psi(\theta)\right \rangle +(r-\varepsilon \left \langle N_\psi (\theta)\hspace{1pt},\hspace{1pt} (\cos 2\pi x, \sin 2\pi x) \right \rangle ) .
\end{equation*}
%\expand{We have \[\dot h =\langle N_\psi(\theta), (\dot x,\dot y)\rangle  = \langle N_\psi(\theta),  \psi (\theta)+rN_\psi (\theta)- \psi (y) -\varepsilon(\cos 2\pi x, \sin 2\pi x)\rangle\] from what we get desired equation. }
The second term is positive as $r>\varepsilon$ and the scalar product of two unit vectors is bounded by one. The first term is non-negative as the angle of the internal normal of an ellipse at a point with a chord containing that point is acute. Thus, $\dot h>0$.  

We can deform $h$ slightly so that $\tilde h$ has rational coefficients in $x$, $y$, while preserving $\dot {\tilde h}>0$. Then $\tilde h = 0$ gives a global cross-section, so the flow has no Reeb component and assumption 4c is satisfied.

\subsection{Flow in \texorpdfstring{$|\Omega_y|<1-\delta$}{Omegay < 1-delta}}
\label{sec:flow_|omega_y|<1-delta}

Next, we examine the flow for parameters in a strip  $|\Omega_y|<1-\delta$ with $\delta = K\varepsilon$ for a large enough $K$ (chosen independently of $\varepsilon$). 

In this strip we can deduce a lot from the case $\varepsilon = 0$, even though its resonance region degenerates to the ellipse $\psi$. The flow has two horizontal invariant circles at the roots of $\sin 2\pi y= \Omega_y$.
%\expand{ We have \[\begin{split}    \dot x &= \Omega_x-\cos2\pi (y-\phi)\\    \dot y &= \Omega_y-\sin 2\pi y\end{split}\]    Now we study if a $y= y^*$ is invariant and we find it is if $\dot y |_{y=y^*} =0$ that is if $\sin 2\pi y^* = \Omega_y$ as stated above (there we don't bother on introducing $y^*$). }
The one with $|y|<\frac{1}{4}$ is linearly attracting; the one with $|y-\frac{1}{2}|<\frac{1}{4}$ is linearly repelling. 
 They have rates
\[\frac{\dd}{\dd y}(\Omega_y-\sin 2\pi y) =\mp 2\pi \sqrt{1-\Omega_y^2},\]
%\expand{As the system on $y$ is decopled of $x$  we can compute ratio convergence (maybe up to a constant depending how exactly you define rate of convergence) as eigenvalue of 1 dim system on $y$. So we get rate of convergence \[\frac{\dd}{\dd y}(\Omega_y-\sin 2\pi y)|_{y=y^\ast}= -\cos 2\pi y^*  =\mp 2\pi \sqrt{1-\Omega_y^2}\] as $\sin 2\pi y= \Omega_y$.  }
which exceed $2\pi \sqrt{\delta}$  in magnitude. 
%\expand{Enough to check inequality for their squares, i.e. $4\pi ^2 (1-\Omega_y^2)>4\pi ^2 \delta$. This is satisfied as    \[1-|\Omega_y|^2>1-(1-\delta)^2 = 2\delta -\delta^2\]   which is bigger than delta as $\delta -\delta^2>0$ for $\delta $ small (it is as K independent of epsilon, so can think first choose K then epsilon small enough). } 
The horizontal velocity is $\dot x =  \Omega_x \mp C\sqrt{ 1- \Omega_y^2}- S\Omega_y,$
%\expand{\[\dot x|_{y= y^*} = \Omega_x-C\cos2\pi y^* -S\sin 2\pi y^*  \] and the substituding $\sin 2\pi y^*= \Omega_y$ leads to desired equation}
so to the right (resp. left) of $\psi$ we have two right-going (resp. left-going) periodic orbits. In the hole the attracting one is left-going and the repelling one is right-going. 
%\expand{Note that the po is where $\dot y =0$, so we will also have $\dot x=0$ there only if we have equilibria, i.e. at resonance (where we actually have a curve of equilibria). So the value of $\dot x$ can only changes when crossing $\psi$ and can check its signs in each region for $\Omega_y=0$ and changing $\Omega_x$ (note that there are 3 regions as we are in the band $|\Omega_y|<1-\delta$). Can alternative argue by taking equational formulation of elipse.   } 

\begin{figure}[htbp]
    \centering
        {\includegraphics[width=0.8\textwidth]{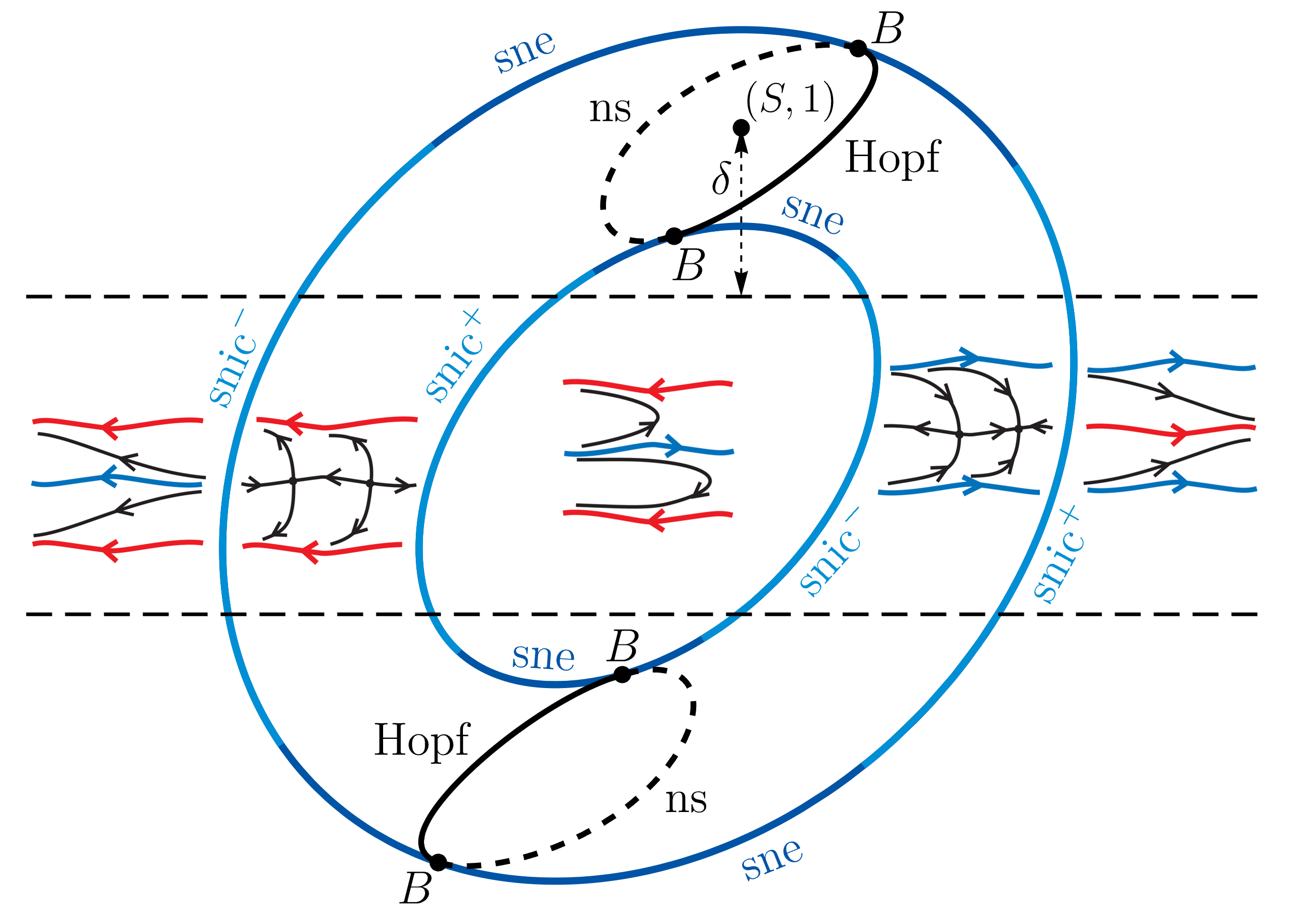}}
    \caption{Phase portraits in the five regions of parameter space into which the $\textnormal{snic}$  curves divide the strip $|\Omega_y|<1-\delta$. Attractive periodic orbits are depicted in red and repelling ones in blue. The signs $\pm$ indicate that the  periodic orbit created by $\textnormal{snic}$ is  right-going/left-going.
    The transitions to $\textnormal{snic}$ along the curves of $\textnormal{sne}$  require $Z$ points that are not depicted here. 
    %Note: Would be more consistant to define sign as what direction the connection goes. However it is not a homoclinic connection (there is no saddle) so tricky to name, and I think for understanding picture is better to think about po they create (so you can match phase diagram outside with snic). 
    %Note: Despite its depiction here this strip covers most of the resonance region for $\varepsilon$ small. 
    }
    %NOTE: using phi = 1.1/24 and eps=0.3
    \label{fig:delta_band}
\end{figure}

By normal hyperbolicity, the $\mathcal C^1$ invariant circles persist with their stability for small $\varepsilon$ if $K$ is chosen large enough.
%\expand{LLL PENDING. Add what exactly needs checking and how the magnitudes are bigger. Need to compute only in unperturbed? or also in perturbed?. Also do you need to use some sort of compactness argument that we get global epsilon for which we have consequence of normal hyperbolicity for all parameters (specially when close to resonance)?}
Moreover, the rest of the flow goes from the repellor to the attractor, as the vector field has non-zero vertical component away from them and they must contain any invariant set for a neighbourhood.
%\expand{In a neighbourhood the invariant circles are the unique invariant sets (together equilibria in them and arcs between them) from nomal hyperbolicity. This is a generic consequence as for hyperbolic equilibria unique in a neighbourhood. Vertical component away from $\sin 2\pi y= \Omega_y$ is non-zero, as there is only place when it zero for the unperturbed.  LLL PENDING: Need to check that we can choose the neighbourhood close and away in such a way that we cover all the phase diagram (for all parameters). Not sure how you quantify where normal hyperbolicity applies.  }
In the hole (restricting to $|\Omega_y|<1-\delta$), these invariant circles are periodic orbits with opposite directions, thus bounding two Reeb components. Outside $\mathcal R$ they are periodic orbits with the same direction, thus forming a Poincar\'e flow, see  figure~\ref{fig:delta_band}. In particular, assumptions 4a and 4b are satisfied in the strip $|\Omega_y|<1-\delta$.

If in $\Omega_x\geq0$ we move our parameters from the hole to the outside of $\mathcal R$, we see that the attracting periodic orbit experiences a saddle-node on an invariant circle (snic) bifurcation on entering the resonance region, then the equilibria do roughly half a revolution in opposite directions, meeting and experiencing another snic bifurcation when leaving the resonance region, see figure~\ref{fig:delta_band}. We have the same behaviour for the repelling circle in $\Omega_x\leq 0$.

To make the transition along the  $\textnormal{sne}$ curves from the $B$ points to the arcs of $\textnormal{snic}$, there must be a $Z$ point, whose study is deferred to the next section.

\subsection{ First look at flow in \texorpdfstring{$|\Omega_y\mp 1|\leq \delta$}{Omegay mp 1 < delta} and \texorpdfstring{$|\Omega_x\mp S|\leq \sqrt{3\delta}$}{Omegax mp S < sqrt(3 delta)}}
\label{sec:Omega_y_mp_1|_leq_delta}
We now consider all parameter values not studied in the previous sections that are close to the resonance region. Thus, we are interested in finding where the resonance region intersects the bands $|\Omega_y\mp 1|\leq \delta$. Using that the resonance region is an $\varepsilon$-neighbourhood of the ellipse $\psi$ we find that
 \[|\Omega_x\mp S|\leq C\sqrt{2\delta(1+1/K)}+O(\varepsilon),\]
so for $\varepsilon$ small enough and $K>2$ we have $|\Omega_x\mp S|\leq \sqrt{3\delta}$. The factor 3 is missing in \cite{Baesens_2018} owing to an oversight. Without loss of generality, from now on we will study the top of the resonance region, i.e. $|\Omega_y- 1|\leq \delta$ and $|\Omega_x- S|\leq \sqrt{3\delta}$, which we denote by $\mathcal Q$. 

Inspired by \eqref{eq:ellipse_parametrisation} we consider the parametrization, 
\begin{equation}
\label{eq:parametrization_sqrt(delta)_negih}
\Omega = \psi (\theta) + \varepsilon \rho N_\psi (\theta),
\end{equation}
so the resonance region is $|\rho|\leq 1$. Writing \[\theta = \frac{1}{4}- \sqrt{\varepsilon}\,\alpha,\]
one finds on Taylor expanding  \eqref{eq:parametrization_sqrt(delta)_negih} that
 \begin{equation}
 \label{eq:Omega_Taylor}
 \begin{split}
     \Omega_x &= S +2\pi \sqrt{\varepsilon} C\alpha  -2 \pi^2 \varepsilon S \alpha^2 + O(\varepsilon^{\frac{3}{2}})\\
 \Omega_y &= 1+ \varepsilon(\rho -2\pi^2\alpha^2) + O(\varepsilon^2)
 \end{split}
 \end{equation}
so in $\mathcal Q$ we have $\rho , \alpha=O(1)$. 
Now as $\Omega_y=1+O(\varepsilon)$, outside a neighbourhood of order $\sqrt{\varepsilon}$ of $y=\frac{1}{4}$, $\dot y >0$, and the flow eventually reaches this neighbourhood.
%\expand{We can be more precise.  Looking at the $\eta$-nullcline below  we find that for \[\eta^2 > \alpha ^2 -\frac{\rho}{2\pi^2} + \frac{\sin 2\pi x}{2\pi^2} + O(\varepsilon)\] we have $\eta'> 0$. Then $\eta$ will only increase and it is also clear that $y$ will increase (even in the region not covered by $\eta$, simply evaluate taylor on $y$ there). So to find the neigh of $y$ it is enough to bound the rhs, which we can do as from  \eqref{eq:Omega_Taylor}  we have that $|\Omega_y -1|\leq \delta $ corresponds to $|\rho -2\pi^2\alpha^2|<K+O(\varepsilon)$ so that the neighbourhood is \[ \eta <\sqrt{\frac{K}{2\pi^2} + \frac{1}{2\pi^2}}+O(\varepsilon) .\] }
%NOT NEEDED BUT USEFUL TO KNOW: \expand{From \eqref{eq:Omega_Taylor} we see that  $|\Omega_x\mp S|\leq \sqrt{3\delta}$  corresponds to $|\alpha|< \frac{\sqrt{3K}}{2\pi C}+O(\sqrt{\varepsilon})$. Then $|\Omega_y -1|\leq \delta $ corresponds to $|\rho -2\pi^2\alpha^2|<K+O(\varepsilon)$ which from the previous bound gives $|\rho|<(1+\frac{3}{2C^2})K +O(\sqrt{\varepsilon})$.}
Thus it is enough to study the flow in the regime $y= \frac{1}{4} + O(\sqrt{\varepsilon}) $ which we do by writing 
\[y=\frac{1}{4}+\sqrt{\varepsilon}\,\eta,\]
and as in this region the motion is slow we take $s=\sqrt{\varepsilon} \,t$ and use $\primebullet{\vphantom{a}}$ for $\frac{\dd}{\dd s}$. Then, Taylor expanding \eqref{eq:principal}  in these coordinates we get, 
\begin{subequations} \label{eq:approx_with_alpha}
    \begin{align}
    \primebullet x  &= 2\pi C(\eta +\alpha) + \sqrt{\varepsilon}\left(2 \pi^2 S ( \eta^2 -\alpha^2) -    \cos 2 \pi x \right ) +O(\varepsilon)\label{eq:approx_with_alpha_x'}\\
        \primebullet\eta&=2\pi ^2(\eta^2-\alpha ^2)+\rho-\sin 2\pi x+O(\varepsilon).
        \label{eq:approx_with_alpha_eta'}
    \end{align}
\end{subequations}
%\expand{We have \[\cos 2\pi y = \cos 2\pi (1/4+\sqrt{\varepsilon}\eta) = -\sin 2\pi \sqrt{\varepsilon}\eta = -2\pi \sqrt{\varepsilon} \eta + O (\varepsilon^{3/2})\]         \[\sin 2\pi y = \sin 2\pi (1/4+\sqrt{\varepsilon}\eta) = \cos 2\pi \sqrt{\varepsilon}\eta = 1-2\pi ^2\varepsilon \eta^2 + O(\varepsilon^2).\]  Substituding this together with the taylor in the parameters we get \[x'= \frac{1}{\sqrt{\varepsilon}}\dot x = \frac{1}{\sqrt{\varepsilon}}(\Omega_x - C\cos2\pi y - S\sin s\pi y - \varepsilon \cos 2\pi x) = \]       \[\frac{1}{\sqrt{\varepsilon}} \left (S + 2\pi \sqrt{\varepsilon}C \alpha - 2\pi ^2\varepsilon S\alpha ^2 - C (2\pi \sqrt{\varepsilon}\eta)- S (1-2\pi^2\varepsilon\eta^2)-\varepsilon \cos 2\pi x+O(\varepsilon^{3/2})\right ) \]    Note that $S$ terms get cancelled and we get desired result. Similarly,     \[\eta' = \frac{1}{\varepsilon}\dot y = \frac{1}{\varepsilon}(\Omega_y-\sin 2\pi y - \varepsilon\cos 2\pi x) = \]        \[\frac{1}{\varepsilon}\left (1 + \varepsilon(\rho - 2\pi ^2\alpha  ^2) - (1-2\pi^2\varepsilon\eta^2) - \varepsilon\sin 2\pi  x+O(\varepsilon^2)\right )\]      the 1's cancel each other and we get desired result. }
    
In \cite[Section 2.7]{Baesens_2018} it was shown that the study of the first order approximation\footnote{In \cite{Baesens_2018} the case $\phi=0$, so $C=1$, was considered but the analysis is no different for $\phi\in (\frac{1}{24},\frac{5}{24})$.} of \eqref{eq:approx_with_alpha}  leads to the following conclusions.

\begin{figure}[htbp]
    \centering
        {\includegraphics[width=0.93\textwidth]{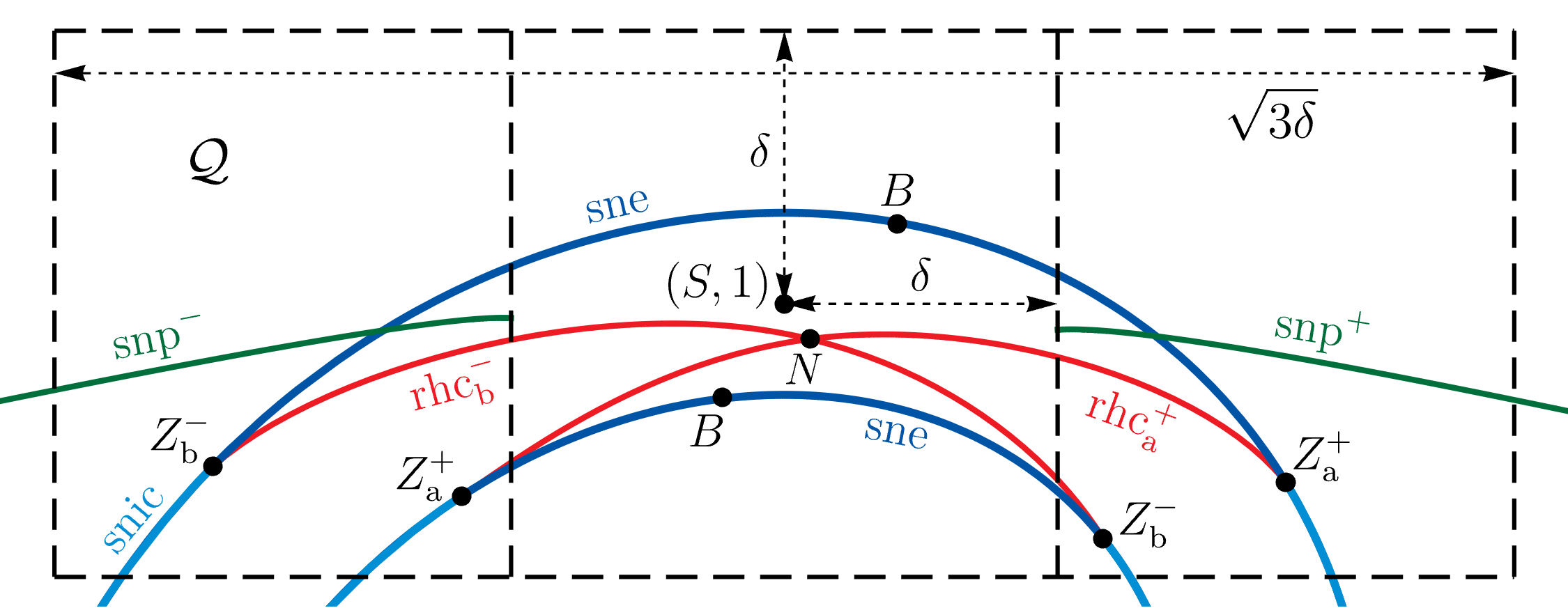}}
    
    \caption[Caption for top resonance]{Sketch of the bifurcations diagram found so far in the top of the resonance region, $\mathcal Q$, excluding the trace-zero curves. We also depict the $\delta$-neighbourhood of $\vect \Omega= (S,1)$ in dashed lines. The sign $\pm$ indicates right-going/left-going and  a/b indicates that the homoclinic connection is above/below the saddle\protect\footnotemark. 
  }
  %NOTE: snp curve should be "increasing" (should eventually get higher than maximum of sne curve, since much left or right in regions 5 and 6 in BM18 we are below them). 
  %EXPAND: sketch is use to illustrate that the figures are not drawn exactly (so ploted as convinient not by formule)
  %NOTE: actually using phi=0.7 to increase readability. Also using smaller epsilon = 0.15 than before.
  %\expand{Note that homoclinoic connection should only be use for saddles as technically cannot be applied to Z, but we still have a orbit connecting equilibria to itself so clear meaning. Also what needs to be right/left-going is the feauter, so the periodic orbit for snp and the homoclinic connection for rhc.  }
    \label{fig:top_resonance}
\end{figure}
\footnotetext{More formally, the angle less than $\pi$ in the rhc at the saddle is above/below the rhc.}

\begin{itemize}[leftmargin=12pt,labelsep=3pt]
    \item For parameter values in the intersection of $\mathcal Q$ and the hole, the flow consists of precisely two Reeb components. Putting this together with the bottom of the resonance and section \ref{sec:flow_|omega_y|<1-delta} we conclude that our example satisfies assumptions 4a and 4b. 

\item In the top of the resonance region there are precisely two curves where the vector field has a generic horizontal rotational homoclinic connection (rhc). Moreover, these curves are given by smooth graphs $\alpha^\pm(\rho)$. There are no other horizontal rhc curves, besides the analogous ones in the bottom of $\mathcal R$ so  assumption 5c is satisfied. The two rhc curves from the top of the resonance region intersect at $\alpha = O(\sqrt{\varepsilon})$ forming an $N$ point, i.e. a saddle with two rhc. Uniqueness of this intersection will be deduced in section \ref{sec:rhc}. 

\item Each horizontal rhc curve meets each boundary curve of $\mathcal R$ once where we have a generic $Z$ point,
i.e.~a saddle-node with a rhc between the neutral direction and the hyperbolic one. So we have precisely\footnote{The outer boundary may have more $Z$ points of other homotopy types, as depicted in figure~\ref{fig:full_bif_diag}. } four $Z$ points of horizontal homotopy type on each component of the boundary of $\mathcal R$,
thus justifying assumption 5\b. It also justifies assumption 5a, as in the inner boundary we have horizontal homotopy type.  
%$ NOTE: I think as rhc generic and sne generic this Z point is also generic but this was not specificly mentioned in BM18

\item Assuming genericity, we can deduce that in the  band $|\rho-2\pi^2\alpha^2|<1$ and to the right (resp. left) of a $\delta$-neighbourhood of $\vect \Omega= (S,1)$ (in the maximum norm),  we have a curve of  horizontal rotational saddle-node periodic orbits (snp) above the corresponding rhc curve, see figure~\ref{fig:top_resonance}. The snp could possibly have cusps and self-intersections, though they would imply regions with more than two horizontal periodic orbits. We did not rule this out but consider it unlikely; moreover it does not directly affect any of the properties we want our system to satisfy. 
\end{itemize}

\subsubsection{Lack of chc and cpo  outside  a $\delta$-neighbourhood of $\vect \Omega= (S,1)$ while in $\mathcal Q$. }
\label{sec:lack_of_chc_cpo}
In the forthcoming section \ref{sec:delta_neigh} we will narrow our attention to a $\delta$-neighbourhood of $\vect \Omega= (S,1)$ in the maximum norm, see figure~\ref{fig:top_resonance}. So here we show that in $\mathcal Q$ but outside a $\delta$-neighbourhood of $\vect \Omega= (S,1)$, there is no parameter value with a contractible homoclinic connection (chc) nor  a contractible periodic orbit (cpo). Using \eqref{eq:Omega_Taylor} we find that in this region $|\alpha|>\sqrt{\varepsilon}\frac{K}{2\pi C}+O(\varepsilon) $ 
%\expand{ We actually have $\frac{K}{2\pi C}\sqrt{\varepsilon}+O(\varepsilon)<|\alpha| <\frac{\sqrt{3K}}{2\pi C}+O(\sqrt{\varepsilon})$. To proof it, from \eqref{eq:Omega_Taylor} we have,     \[\Omega_x-S = 2\pi C \alpha (1-\pi\frac{S}{C}\alpha \sqrt{\varepsilon}  )\sqrt{\varepsilon} + O(\varepsilon^{3/2}).\]     So as the left hand side is at least $\delta = K\varepsilon$ we find        \[\frac{K}{2\pi C}\sqrt{\varepsilon} + O(\varepsilon)<  |\alpha| (1-\pi\frac{S}{C}\alpha \sqrt{\varepsilon}  ) .\]       (Note that as we have error term we can use $<$).   Dividing by the parenthesis while using that $\alpha = O(1)$ and  $\frac{1}{1-r}= 1+r +O(r^2)$ we find that $\frac{K}{2\pi C}\sqrt{\varepsilon} + O(\varepsilon)<\alpha$. For the upper-bound one less order is needed so that \[\Omega_x-S = 2\pi C \alpha \sqrt{\varepsilon} + O(\varepsilon),\] and as the lhs is at most $\sqrt{3\delta}= \sqrt{3K}\sqrt{\varepsilon}$ we get the desired result. }
and as cpo and chc require equilibria, we can restrict our attention to $|\rho|\leq 1$. 
%\expand{as equilibria only exist in the resonance region which is given by this. }

First note that the divergence of \eqref{eq:approx_with_alpha} (i.e.~trace of the derivative of the vector field) is
\begin{equation}
    \tr = 4\pi^2\eta + 2\pi \sqrt{\varepsilon}\sin 2\pi x+ O(\varepsilon).
    \label{eq:tr}
\end{equation}
From  \eqref{eq:approx_with_alpha_x'} we find that the equilibria are located at $\eta = -\alpha + \sqrt{\varepsilon}\frac{\cos 2\pi x}{2\pi C}+O(\varepsilon)$ 
%\expand{. From $x'=0$ is clear that $\eta = -\alpha+O(\sqrt{\varepsilon})$. So let $\eta = -\alpha+a\sqrt{\varepsilon}+O(\varepsilon)$. Substituting back in $x'=0$ we get        \[0=2\pi C a\sqrt{\varepsilon} + (2\pi^2 S(-2\alpha a\sqrt{\varepsilon}+O(\varepsilon))-\cos 2\pi x)\sqrt{\varepsilon}+O(\varepsilon)\]        and thus we get $a= \frac{\cos 2\pi x}{2\pi C}$.}
so that the trace at them is 
\[\tr = -4\pi^2\alpha + 2\pi \sqrt{\varepsilon}\left (\sin 2\pi x +C^{-1}\cos 2\pi x \right )+ O(\varepsilon)\]
The sign of this expression is then dominated by the first term if $K>2(C+1)$ as $|\alpha|>\sqrt{\varepsilon}\frac{K}{4\pi C}$ for $\varepsilon$ small enough. 
%\expand{From $|\alpha|>\frac{K}{2\pi C}\sqrt{\varepsilon}+O(\varepsilon) $ add factor 2 to get rid of O term. Then we need       \[4\pi^2\alpha>  2\pi \sqrt{\varepsilon}|\sin 2\pi x +C^{-1}\cos 2\pi x |\]        is enough to check      \[4\pi^2\frac{K}{4\pi C}\sqrt{\varepsilon}>  2\pi \sqrt{\varepsilon}|\sin 2\pi x +C^{-1}\cos 2\pi x |\]         simplifying and bounding absolute value by $(1+C)$ we get       \[\frac{K}{ C} > 2 (1+C^{-1})\]     which leads to the desired result.     }

Now a cpo, $\gamma$, with period $T$ would have Lyapunov exponent $\frac{1}{T}\int_\gamma \tr \dd s$ which from \eqref{eq:tr} and  \eqref{eq:approx_with_alpha_x'} can be written as,
\begin{equation}
    -4\pi^2 \alpha +\frac{2\pi}{C} \sqrt{\varepsilon} \int_\gamma\frac{1}{T} \biggl ( \frac{\primebullet x}{\sqrt{\varepsilon}}  - 2\pi^2 S ( \eta^2 - \alpha^2) + \cos 2\pi x + C \sin 2\pi x \biggr )\,\dd s+ O(\varepsilon).
    %NOTE: not use left right as they leave more space in parenthesis and then eq number doesnt fit
    \label{eq:lyapunov_cpo}
\end{equation} 
%\expand{Substituding $\tr$ we get      \[\frac{1}{T}\int_\gamma \tr \dd s =\frac{1 }{T}\int_\gamma 4\pi^2\eta + 2\pi \sqrt{\varepsilon}\sin 2\pi x+ O(\varepsilon)\dd s\]        Then substituding $\eta$ in $x'$ equation (while using the fact that we are normalizing by time, so taking average, so $O$ terms and constants leave the integral) we get the result.     }
Again we want to show that the sign of the Lyapunov exponent is dominated by the first term, so we focus on bounding the integral.  To this end  we note that the $\primebullet\eta$-nullclines are contained in  
\begin{equation}
    2\pi^2|\eta ^2-\alpha^2|< |\rho|+1+O(\varepsilon),
    \label{eq:eta_nullcline_bound}
\end{equation}
%\expand{$\eta$-nullcline  are at     \[2\pi^2(\eta ^2-\alpha^2)= \sin  2\pi x + \rho+O(\varepsilon),\]     So modulus lhs is same as rhs which is bounded by $1+|rho| + O(\varepsilon)$. Note: We CANT do it without finding maximum and minimum rhs, as for the argument we do later we NEED horizontal bands. Note: From null clines we could isolate $\eta^2$ and then we get a bound on it in terms of order $\sqrt{K}$ (by substituding upper bound on $\alpha$). This in fact quantifies the $y$ neighbourhood of order $\sqrt{\varepsilon}$ that we discuss in the start of this section.}
so outside these horizontal bands $\eta$ is strictly monotonic and thus  cpos must be contained in \eqref{eq:eta_nullcline_bound}.
%\expand{If a cpo would intersect region outside, it would be in a horizontal region with $\eta $ monotonic which is impossible as orbit must leave the horizontal region (otherwise orbit monotone) but then when comes back to the region $\eta$ must have opposite monotonicity as the boundary is horizontal which contradicts monotonicity eta in that region.   (could have rotational po doing lap) }
Moreover, as $\gamma$ is contractible,  $\int_\gamma\primebullet x \,\dd s=0$ so that the integral in \eqref{eq:lyapunov_cpo} is bounded above by 
\[\frac{1}{T}\int_\gamma S(|\rho|+1) + |\cos 2\pi x |+C|\sin 2\pi x| \,\dd s + O(\varepsilon), \] 
and hence also by  $2S + C +1$. 
%\expand{Seems that we are ignoring $O(\varepsilon)$ term, but we are not as we are approximating integral $\sin$ by all 1, which gives enough leeway. }
 Recall that $|\alpha|>\sqrt{\varepsilon}\frac{K}{4\pi C}$ so that if $K>2(2S+C+1)$, the sign of the Lyapunov exponent is given by the first term in \eqref{eq:lyapunov_cpo}. 
%\expand{Plugging computation of integral and $\alpha$ in \eqref{eq:lyapunov_cpo} it is enough to show that       \[4 \pi ^2\frac{K}{4\pi C} \sqrt{\varepsilon}=   \frac{\pi K}{ C} \sqrt{\varepsilon}>\frac{2\pi}{C}(2S + C +1)\sqrt{\varepsilon}\]      (Can ignore $O$ terms as aproximation we took from $\alpha$ wasn't tight.) From this we get the result.  }

In conclusion, the sign of the Lyapunov exponent coincides with the sign of the trace at the equilibria.  We can deduce then that there cannot exist cpos in this regime, as they would have the same attracting/repelling nature as the node and would give rise to impossible phase portraits. Essentially the same argument shows that there are no chc either.
%EXPAND: cpo needs to have node inside (no saddle or saddlenode by index) then if both attracting need cpo (as no other equilibria)in between repelling (but it would be attracting). For chc note that as saddle generic also need node inside, so same argument (saddle node would also need node inside as it is not type of degenerate equilibria that allows exclussibly petals).

%OBS: Note that standard argument of bendixond dulac doesnt work as actually sign of trance changes inside possible po
%OBS: we technically not care if there is chc outside delta neigh as it would never be able to intersect with ns which are in delta neigh. Still it is nice to know that there aren't any.  

\section{\texorpdfstring{$\delta$}{delta}-neighbourhood of top and bottom of resonance region}
\label{sec:delta_neigh}

To complete our analysis of the vector field, we zoom in on the $\delta$-neighbourhood of $\vect \Omega=(S,1)$ in the maximum norm as depicted in figure~\ref{fig:top_resonance}(similar analysis can be done for $\vect \Omega=(-S,-1)$).
Here, we derive an approximation for the vector field which is reversible (conjugate to 
its time-reverse) and symplectic (locally Hamiltonian).  We use the remainder terms to prove that for small $\varepsilon$ there is a unique 
$N$ point which is connected to the $B$ points by arcs of chc. We will also use these remainder terms in the upcoming section \ref{sec:cpo} and \ref{sec:K_H_points} to show that there are precisely two $K$ points and a 
unique $H$ point, and give numerical evidence that there is at most one cpo.

\subsection{Reversible approximation}
\label{sec:reversible_approx}
 To study a $\delta$-neighbourhood of $\Omega=(S,1)$ we specialize the approximation in \eqref{eq:approx_with_alpha} by letting $ \alpha = \sqrt{\varepsilon}\,\tilde \alpha /(2\pi)$, 
% \expand{Note that then $\Omega_x - S = \varepsilon C \tilde \alpha +O(\varepsilon^{3/2})$ and $\Omega_y-1= \varepsilon \rho +O(\varepsilon^2)$, so that $\delta$ neighbourhood is given by $|\rho|<K+O(\varepsilon)$ and $|\tilde \alpha |< K/C + O(\sqrt \varepsilon)$  }
 so
 \begin{equation}
 \label{eq:reversible_approx_sqrt_eps}
\begin{split}
    \primebullet x &= 2\pi C \eta + \sqrt{\varepsilon}\,(C\tilde \alpha- \cos 2\pi x + 2 \pi^2 S\eta^2 ) + O(\varepsilon)\\
    \primebullet \eta &= \rho -\sin 2\pi x + 2 \pi^2 \eta^2 + O(\varepsilon).
\end{split}
\end{equation}

In this section we study the unperturbed case $\varepsilon=0$,
\begin{equation}
\label{eq:reversible_approx}
\begin{split}
    \primebullet x &= 2\pi C \eta \\
    \primebullet \eta &= \rho -\sin 2\pi x + 2 \pi^2 \eta^2 
\end{split}
\end{equation}
which is a family of reversible vector fields with respect to the reflection $\eta \rightarrow -\eta$ and we denote it by $\mathbf v_\rho$. Inspired by \cite[Section 2.8.1]{Baesens_2018}, we find that on the universal cover they are Hamiltonian with respect to the symplectic form,
\[\mathcal A=e^{-2\pi x/C}  \dd x\wedge \dd\eta .\]
%\expand{The actual way to find this factor is to consider the change in standard volue which is given by trace so $V ' = 4\pi^2 \eta V $ then substituding $x'$ get $V ' = 2\pi/C x' V $ so $\frac{d V}{dx} =2\pi/C V $ with solution $V = V_0 e^{2\pi/C x }$ so that volume increases at this rate and thus to be hamiltonian (conserves volume) need to compensate this by mutiplying by inverse.}
%\expand{It is simplectic as is not degenerate and obviously it is closed as there is nothing in dimension 3. }
Indeed, we have
\[e^{-2\pi x /C}\,\primebullet x = H_\eta, \hspace{1.6cm} e^{-2\pi x /C}\,\primebullet \eta = -H_x,\]
where the Hamiltonian is
\begin{equation*}
H(x,\eta)=e^{-2\pi x /C}\left( \pi C \eta^2+\frac{C\rho}{2\pi} - \frac{1}{4\pi}g(x) \right),
\label{eq:Ham}
\end{equation*}
%\expand{From integrating  $e^{-2\pi x /C}x' = H_\eta$ get $H = e^{-2\pi x/C}(\pi C \eta^2+\tilde g)$. Then from  $e^{-2\pi x /C}\eta' = -H_x$ we get $\frac{2\pi }{C} \tilde g(x)-\tilde g'(x)=\rho -\sin 2\pi x$. Taking $\tilde g(x) = \frac{C}{2\pi}\rho - \frac{1}{4\pi }g(x)$  we get      \[\frac{1}{2C}g(x)-\frac{1}{4\pi}g'(x)=\sin 2 \pi x\]    One can solve this with inegration factor method, but easier to find particular solution. Trying $g(x)= A\cos 2\pi x+B\sin 2\pi x$ or equation becomes \[A-CB=0 \hspace{1cm} B+CA\] from what we attein the desired result. Note that we could also add to this solution an homogeneous one $\tilde C e^{2\pi x/C}$ when plugging this into $H$ it is equivalent to adding a constant to $H$ as this term is inside the parenthesis.   }
with 
\begin{equation}
g(x)=\frac{2 C}{1+C^2}\left( C \cos 2\pi x + \sin 2\pi x\right ),
\label{eq:def_g}
\end{equation}
which is found by solving the relation $\frac{1}{2C}g(x)-\frac{1}{4\pi}g'(x)=\sin 2 \pi x$.  Because $C>0$, this can be rewritten as 
\[g(x)=\frac{2 C}{\sqrt{1+C^2}}\cos 2 \pi (x -x_g)\]
%\expand{We have \[r\cos 2\pi (x-x_g) = r\cos 2\pi x_g \cos 2 \pi x +r\sin 2\pi x_g \sin 2 \pi x \] so that $r$ and $x_g$ are sort of the polar coordinates to $A$ $B$ from the previous expand and we have $r= \sqrt{A^2 + B^2} $ and $2\pi x_g = \arctan B/A$ which leads to the desired result. Note that it is standard determination as the coefficients $A$ and $B$ are positive.  }
with $x_g=\frac{1}{2\pi}\arctan \frac{1}{C}$, where we take the standard determination of arctangent from $-\frac{\pi}{2}$ to $\frac{\pi}{2}$.

\begin{figure}[htbp]
    \centering
        {\includegraphics[width=0.78\textwidth]{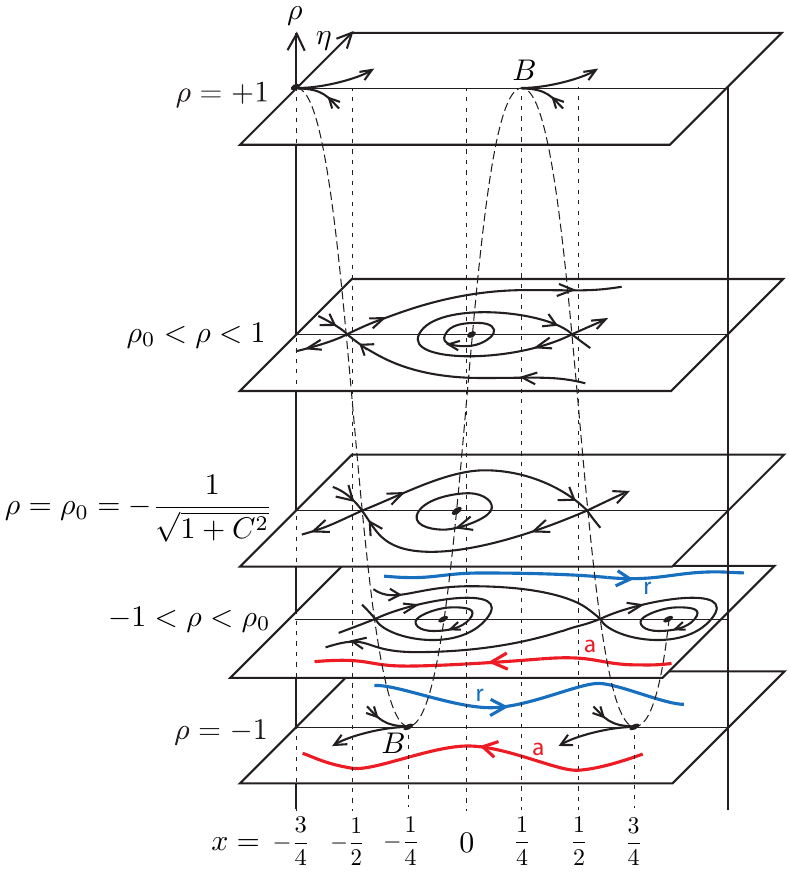}}
    \caption{Sketch of the phase portraits for the reversible family of vector fields \eqref{eq:reversible_approx} at selected 
values of $\rho$. To increase visibility, the $\rho$-axis and the phase portraits have been slightly distorted.}
%EXPAND: sketch is use to illustrate that the figures are not drawn exactly (so ploted as convinient not by formule)
    \label{fig:hamiltonian}
\end{figure}

We now study the phase portraits using the fact that $H$ is conserved along trajectories. They have two distinct regimes separated by a parameter value $\rho_0$ with a necklace
%EXPAND: actually Neckless is a parameter value so we are being a bit loose here 
as depicted in figure~\ref{fig:hamiltonian}. To compute $\rho_0$ we first note that as our phase space is a cylinder $\mathbb T^1\times \mathbb R$ any rotational homoclinic connection or periodic orbit will be invariant under $x\mapsto x+1$.  
%\expand{They have to come back same point after some revolutions, so in principle could have $x->x+n$ with $n$ integer. After one revolution either we connect back to where we where or we are above/below, in which case we will always be above/below and can't connect (consequence of being specifically on cylinder). }
As $H $ has a non-periodic factor, the only invariant level set under $x\mapsto x+1$ is $H=0$. 
%\expand{Must have energy 0 otherwise change in $x$ will afect it. Moreover, as all the rest is periodic, it will indeed be invariant. }

The equilibrium points of the system are at $\eta =0$ and
$\sin 2\pi x=\rho$, where $\det = 4\pi ^2 \cos 2\pi x$, so the saddle is at $x_{\sad} \in (\frac{1}{4},\frac{3}{4})$ mod 1 and the centre at $x_{\cen} \in (-\frac{1}{4},\frac{1}{4})$ mod 1. The rhc occurs when the energy at the saddle,  $E_{\sad}=H(x_{\sad},0)$, vanishes, i.e. 
\[\frac{C\rho_0 }{2\pi} = \frac{g(x_{\sad})}{4\pi}\quad \Longleftrightarrow \quad \rho_0 =\frac{1}{1+C^2}\left (-C\sqrt{1-\rho_0^2}+\rho_0\right )\]
%\expand{Substituded $\sin 2\pi x_{\sad}=\rho$ and $\cos 2\pi x_{\sad} = -\sqrt{1-\rho_0^2}$, negative because in saddle trace (which is also cos) negative. } 
and with some algebra we get $\rho_0=-1/\sqrt{1+C^2}$. 
%\expand{\[\iff C\rho_0=-\sqrt{1-\rho_0^2}\iff C^2\rho_0^2 = 1-\rho_0^2 \iff \rho_0 = -1/\sqrt{1+C^2}\] where we take negative root as $\rho_0$ must be negative by the first equality in this expression. }
The rhc are then given by the level set  $H=0$ at $\rho_0$ so 
\begin{equation}
    \eta = \pm \sqrt{\frac{\cos 2\pi (x-x_g) +1}{2\pi^2\sqrt{1+C^2}} } =\pm \frac{1}{\pi\sqrt[4]{1+C^2}}\cos  \pi ( x -x_g). 
    \label{eq:necklace_unperturbed}
\end{equation}
%\expand{$H=0$ is equivalent to      \[C\pi \eta^2 =- \frac{C\rho_0}{2\pi}+\frac{g(x)}{4\pi}= \frac{C}{2\pi \sqrt{1+C^2}}+ \frac{2C\cos 2\pi(x-x_g)}{4\pi\sqrt{1+C^2}}\]      Taking factor this lead to the first expression in \eqref{eq:necklace_unperturbed}. Then use the typical formula to integrate powers of trigonometric functions $\cos^2\varphi = \frac{\cos 2\varphi +1 }{2}$ in reverse to get desired result.}
For $\rho \in (-1,\rho_0)$, $E_{\sad}<0$, so its level set includes a contractible homoclinic loop to the right of the saddle, whereas for 
$\rho \in (\rho_0,1)$, $E_{\sad}>0$, so it has a contractible homoclinic loop to the left. 
%\expand{The level set is \[C\pi\eta^2 = E_{\sad}e^{2\pi x/C}-\frac{C\rho}{2\pi}+\frac{g(x)}{4\pi}\] plotting the rhs with positive energy we find get sinusoidal curve that generally increases (looked from a far seems exponential). for $x$ small this function is positive in islands corresponding to cpo. Then at some point a local minimum happnes at height 0 which corresponds to the chc. To the left of it the function is always positive and corresponds to the other separatrix. \\ Observation: $H$ is not periodic so $H(x_{\sad},0)$ is not well defined (that is why we don't define $E_{\sad}$ here. However, it is well defined if it is positive or negative (doesn't depend on non-periodic factor).   }
Each contractible homoclinic loop bounds a region foliated by cpos around a centre, see figure~\ref{fig:hamiltonian}.

Note that, despite its Hamiltonian nature on the universal cover, the vector field can have attracting and repelling rotational periodic orbits.
%\expand{Note that if it would be Hamiltonian on the cylinder (not universal cover) this wouldn't happens as it would preserve area. }
Indeed, the level set $H=0$ gives explicitly
\begin{equation}
    2\pi \eta = \pm \sqrt{g(x)/C -2\rho }
    \label{eq:rpo_reversible_approx}
\end{equation}
which for $\rho<\rho_0$ gives two rotational periodic orbits, see figure \ref{fig:hamiltonian}.
%\expand{$H=0$ gives $\pi C \eta^2 = \frac{g(x)}{4\pi}-\frac{C\rho}{2\pi}$ from what we get the equation. This analysis doesn't need equilibria so we still get it for $\rho<-1$. }
Integrating the divergence of \eqref{eq:reversible_approx} by changing integration variable to $x$ we find that they have  Floquet multipliers $e^{\pm 2\pi/C}$ respectively.  
%\expand{Use lyapunov mutiplyers instead of coef as they dont depend on period, they are given by $exp(\int_\gamma \tr \dd s)$. In this case we have \[\int_\gamma \tr \dd s = \int_\gamma 4\pi ^2\eta ^2 =  \int_\gamma 2\pi /C x' \dd s =\pm \int_{x_0}^{x_0+1} 2\pi /C\dd x =\pm 2\pi /C\]    where we have used that the orbit with $\eta$ positive goes from $x_0$ to $x_0+1$ whereas the negative goes in the opposite direction. }
The square of \eqref{eq:rpo_reversible_approx} also gives a periodic orbit for $\rho \in (\rho_0,1)$ but it is just one of a continuum of contractible and neutrally stable ones inside the chc. 
%\expand{Neutrally stable as contractible so in integration in previous expand would be from $x_0$ to $x_0$. (can also think that path with positive $\eta$ is canceled by the negative one). }
% lll need sign value Esad E cen? Answer: I think no direclty although maybe indirectly (for instance in argument showing minimum g is in saddle)

In the following sections we use Pontryagin’s energy balance method, see \cite{Pontryagin_1934} and \cite[Section 2.8.2]{Baesens_2018}, to determine where various features of the reversible approximation persist on adding the correction terms. In this method we consider $\gamma $ a closed orbit of the unperturbed system $\mathbf v_\rho$. Then, the increment of $H$ for nearby orbits in the perturbed vector field $\mathbf v_\rho + \sqrt{\varepsilon} \tilde {\mathbf v} + O(\varepsilon)$ is,
\begin{equation}
    \Delta H = \int_{\gamma} \primebullet H \,\dd s+O(\varepsilon) = \sqrt{\varepsilon}\left ( \int_{\gamma} e^{-2\pi x/C} \tilde v_2\, \dd x -\int_{\gamma} e^{-2\pi x/C} \tilde v_1 \,\dd \eta \right)+O(\varepsilon).
    \label{eq:Delta_H}
\end{equation}
%\expand{\[ \Delta H = \int_{\gamma_\epsilon}  H' \,\dd s= \int_{\gamma}  H' \,\dd s+O(\varepsilon)\]    by Pontryagin. Then \[H'=\langle \nabla H , \mathbf v_\rho + \sqrt{\varepsilon} \tilde {\mathbf v} + O(\varepsilon) \rangle  =e^{-2\pi x/C}(\langle (-v_{\rho 2},v_{\rho 1}), \mathbf v_\rho \rangle + \sqrt{\varepsilon}\langle (-v_{\rho 2},v_{\rho 1}), \tilde {\mathbf v} \rangle ) +O(\varepsilon) \]    \[= \sqrt{\varepsilon} e^{-2\pi x/C} \langle (-v_{\rho 2},v_{\rho 1}), \tilde {\mathbf v} \rangle  +O(\varepsilon) \]     Now as we are integrating over an orbit of the system $(x',\eta')=\mathbf v_\rho$  (its standard velocity, not hamiltonian) so $dx = x'ds = v_{\rho 1} ds$ and $d\eta = \eta'ds = v_{\rho 2} ds$ from what we get desired result. }
so that the closed orbit persists only  when $\Delta H =0$.
Moreover, if the coefficient of $\sqrt{\varepsilon}$ in \eqref{eq:Delta_H} has a non-degenerate zero (or submanifold of non-degenerate zeroes) as a function of parameters, then the orbit only survives along a submanifold in parameter space $\mathcal C^1$-close to where the first term vanishes.
%\expand{Note that no need to worry about multivalued in general for this method.  We are studying persistance of features given by level curves of $H$. So value of $H $ in it doesn't change and we don't get multivalued conection. After perturbation if there is a connection with non-zero increment of $H$ it will be far away from this, so we are not studying it.  \\ Direct argument for our case: in our case if level curve is rotational by definition needs to be at 0 (where $H$ is not multivalued), if is contractible mutlivalued is no relevant.   }

\subsection{Rotational homoclinic connection (\texorpdfstring{$\textnormal{rhc}$}{rhc}) and necklace point (\texorpdfstring{$N$}{N})}
\label{sec:rhc}

First, we compute where the unperturbed $\textrm{rhc}^\pm$ given by \eqref{eq:necklace_unperturbed} persist, by letting $\rho =\rho_0+\sqrt{\varepsilon}\,\tilde \rho$. So the vector field is $\mathbf v_{\rho_0}+\sqrt{\varepsilon}\tilde {\mathbf v} +O(\varepsilon)$ with %$\tilde {v}_2 =\tilde \rho $ 
\[\tilde {\mathbf v} = (C\tilde \alpha- \cos 2\pi x + 2 \pi^2 S\eta^2, \tilde \rho ).\]
%\expand{Note that it is not quite what we explained above with the change in $\rho $ but it is esentially the same. }
If we denote by $\gamma_\pm$  the unperturbed $\textrm{rhc}^\pm$ from $x_g-1/2$ to $x_g+1/2$, we see from \eqref{eq:Delta_H} that they persist along 
\[\Delta H= \sqrt{\varepsilon} \,(c\tilde \rho - Ca\tilde \alpha +b)+O(\varepsilon)=0\]
where,
{\allowdisplaybreaks
 \begin{align*}
   a&=\int_{\gamma_+} e^{-2\pi x/C}  \,\dd \eta=-\int_{\gamma_+} e^{-2\pi x/C} \frac{\sin \pi(x-x_g) }{\sqrt[4]{1+C^2}}\,\dd x =\frac{4C e^{-C^{-1}\arctan{\frac{1}{C}}} \cosh{\frac{\pi}{C}}}{\pi \sqrt[4]{1+C^2}\,(4+C^2)}, \\
    b&=\int_{\gamma_+} e^{-2\pi x/C} ( \cos 2\pi x - 2 \pi^2 S\eta^2) \,\dd \eta\\
    &=-\int_{\gamma_+} e^{-2\pi x/C} \left (\cos 2\pi x-2S\frac{\cos ^2\pi (x-x_g)}{\sqrt{1+C^2}}\right ) \frac{\sin \pi(x-x_g)}{\sqrt[4]{1+C^2}}\dd x\\
    &=-a\frac{2 C^2 (C + 2 S)}{\sqrt{1+C^2} \, (4 + 9 C^2)},\\
    c&=\int_{\gamma_+} e^{-2\pi x/C} \,\dd x=a\frac{1}{4}\sqrt[4]{1+C^2}(4+C^2)\,\tanh \frac{\pi}{C}>0.
\end{align*}
}
%\expand{For $a$ simply do change in variable $\dd\eta = \frac{\dd \eta}{\dd x}\dd x $ where we compute the derivative from \eqref{eq:necklace_unperturbed} so \[-\frac{\sin \pi (x-x_g)}{\sqrt[4]{1+C^2}}\]. For $\gamma_-$ we would take the negative root so would change sign (nothing else would change as the expression inside doesn't depend on eta) justifying the $\pm$ in the expression  $\tilde \rho$. For $b$ first note the sing chosen in the expression of $\tilde \rho$ makes us flip the order of the $\cos$ with the square. Then do same change of variables (again for $\gamma_-$ would simply have negative of this, as the expression only depends on $\eta^2$) and to substitute $\eta^2$ use the first equality in  \eqref{eq:necklace_unperturbed}. To do actual computations we use computer.}
%NOTE: Seems weird to use a instead of c as constant, but we do it because then "complexity" of expressions is more or less the same (the more complex b has a term which is basically a inside integral). Also it is more combinient to do the actual computations we need, as we only need the exact value in N point where we devide by a. 
For $\phi =0$ we recover the coefficients computed in the corrigendum \cite{Baesens_2022_Corrigendum}. Then, by Pontryagin’s method the curves of $\textrm{rhc}^\pm$ are given to leading order in $\sqrt{\varepsilon}$ by 
\begin{equation}
    \tilde \rho =\pm (Ca\tilde \alpha -b)/c,
    \label{eq:rhc_rho}
\end{equation}
and are in fact $\mathcal C^1$-close to \eqref{eq:rhc_rho}. 
 The approximate curves $\textrm{rhc}^\pm$ cross transversely at $\tilde \rho =0$ and 
\begin{equation}
    \tilde \alpha = \frac{b}{Ca}=-\frac{2C(C+2S)}{\sqrt{1+C^2}\,(4+9C^2)},
    \label{eq:tilde_alpha_N}
\end{equation}
forming a generic $N$ point (non-vanishing trace will be shown below). 
%\expand{Crossing is at $\tilde \rho =0$ is only way that plus minus is the same. Then the equations for the lines gives $\tilde \alpha$. They give generic $N$ point as the intersection is transverse  and the trace in saddle is non-zero (we will check we are inside trace zero loop) so that the rhc break in generic manner.}
%Also note for future reference that they cross $\tilde \alpha =0$ at $\tilde \rho=\mp b/c$.

There is no other horizontal rhc in $|\rho|\leq 1$, $\tilde \alpha =O(1)$, because the phase portraits for the reversible approximation allow such a feature only near $\rho =\rho_0$, where Pontryagin’s method above gives the locally unique curves found.
%\expand{In the reversible approax the seperatrix either form chc, rhc, tend to rpo or go to infinity. The only small perturbation of them that can lead to a rhc is rhc or chc/rpo that go very close to the saddle. That happens precisely close to  $\rho =\rho_0$. For this parameters pontryagin detects when we get conenection, if conection is missed in the first lap it can not reconect in second due to top cyclinder (there is above and below).  }
So the $N$ point found above is unique in this region of parameters, which by section \ref{sec:Omega_y_mp_1|_leq_delta} is the only possible location of an $N$ point in the top of the resonance region. Doing the same for the bottom we prove assumption 6.

In addition to the two curves of rhc the unfolding of a necklace point contains two arcs of chc whose direction is determined by the sign of the trace at the saddle  \cite{Baesens_1991}.  
%\expand{By direction we mean direction of the parametric curves (right of left of the N point) but probably ok to think of it as direction of the actual conection as this will determine where the parametric curve can emanete from }
 From \eqref{eq:trace_zero_elipse} and \eqref{eq:Omega_Taylor} we deduce that the trace-zero loop in $(\tilde \alpha, \rho)$ parameters is $\mathcal C^1$-close to the ellipse
%NOTE: C^1-close is to first order 
\begin{equation}
    \rho ^2 + C^2(\tilde \alpha -\rho)^2 = 1.
    \label{eq:trace_zero_in_alpha_rho}
\end{equation}
%\expand{By $\mathcal C^1$-close, we mean that equality is true up to some error $O(\sqrt{\varepsilon})$ and moreover the differential of the error also has this order.  Substituting \eqref{eq:Omega_Taylor}   in  \eqref{eq:trace_zero_elipse}  we get        \[ \left ((S +2\pi \sqrt{\varepsilon} C\alpha  -2 \pi^2 \varepsilon S \alpha^2)- C( 1+ \varepsilon(\rho -2\pi^2\alpha^2)) + (C-S) +O(\varepsilon^{3/2})\right )^2+\] \[+\left ( ( 1+ \varepsilon(\rho -2\pi^2\alpha^2))- 1+O(\varepsilon^{2})\right )^2 = \varepsilon^2.\]     Note that the 0th order term vanishes. Also substituting  $\alpha = \sqrt{\varepsilon}\tilde \alpha /(2\pi)$     \[ \left (  \varepsilon C\tilde \alpha  - C\varepsilon \rho  +O(\varepsilon^{3/2}) \right)^2+\left (  \varepsilon\rho+ O(\varepsilon^{2}) \right )^2 = \varepsilon^2.\]     So to first order we get \eqref{eq:trace_zero_in_alpha_rho} (note all errors come from taylor so $\mathcal C^1$). }
Thus substituting \eqref{eq:tilde_alpha_N} and $\rho_0$ we see that the $N$ point is inside the trace-zero loop.
%\expand{Inside means $\rho ^2 + C^2(\tilde \alpha -\rho)^2 < 1$ (without using higher order terms). We have \[\rho_0 + C^2(\tilde \alpha -\rho_0)^2 =\frac{1}{1+C^2}  + C^2\left (-\frac{2C(C+2S)}{\sqrt{1+C^2}\,(4+9C^2)} +\frac{1}{\sqrt{1+C^2}} \right)^2<1 \]     \[\iff  1+ C^2\left (1-\frac{2C(C+2S)}{(4+9C^2)}\right )^2 <1+C^2\] Now, note that the fraction is clearly positive and smaller than 1, so that it is clear that the square is smaller than 1 and the inequality is satisfied.  }
In particular it is to the right of the ns curve, so the trace at the saddle is negative. Hence, the chc emanate to the right of the necklace point as depicted in   figure~\ref{fig:delta_neighbourhood}.
Note that, as mentioned in the corrigendum \cite{Baesens_2022_Corrigendum}, in the case $\phi=0$ we also get the $N$ point inside the trace-zero loop.

\subsection{Contractible homoclinic connection (\texorpdfstring{$\textnormal{chc}$}{chc})}
\label{sec:chc}
Next we compute how the contractible homoclinic connections of the reversible approximation, $\gamma$, persist for each $\rho \in (-1,\rho_0)\cup (\rho_0,1)$. We apply Pontryagin's method, where the coefficient in \eqref{eq:Delta_H} vanishes if 
\[Ca(\rho)\tilde \alpha -b(\rho) =0 \]
where 
\[ a(\rho)=\int_\gamma e^{-2\pi x/C}\,\dd\eta, \hspace{1.5cm}
    b(\rho)=\int_\gamma e^{-2\pi x/C}(\cos 2\pi x - 2\pi^2 S\eta ^2)\,\dd\eta.
\]
%\expand{In this case we don't zoom in in any $\rho$ so \[\tilde {\mathbf v} = (C\tilde \alpha- \cos 2\pi x + 2 \pi^2 S\eta^2, 0 ).\] and thus we do not get the $dx$ integral. }
So for each $\rho$ there is a unique $\tilde \alpha$ with chc, close to  $\tilde \alpha (\rho) = C^{-1}b(\rho)/a(\rho)$,  since $a(\rho)>0$, as we proceed to show. First note that if we change variables to $x$ and then integrate by parts while using that $\gamma$ is contractible  we get 
\[\int_\gamma e^{-2\pi x/C} h(x)\,\dd\eta = \int_\gamma \left (\tfrac{2\pi}{C}h(x)-h'(x)\right ) e^{-2\pi x/C}  \eta\,\dd x \]
%\expand{First substituting $\dd \eta$ for $\frac{\dd \eta }{\dd x} \dd x$, then $\frac{\dd \eta }{\dd x}$ is the term we integrate rest we differentiate in integration by parts. As $\gamma $ is contractible starting and en point is same saddle (in the universal cover) so "first term" in integration by part vanishes. We could work with normal integrals on x by put initial and end point to be intersection with axes and do the integral in two terms (if we don't want to think about integration on curves). In this case non-integral term is still 0 as we have $\eta=0$.   }
for any function $h$. Picking $h=1$, we obtain, 
\begin{equation}
\label{eq:a_chc}
    a(\rho)= \frac{2\pi}{C}\int_\gamma e^{-2\pi x/C}\eta \,\dd x= \frac{2\pi}{C}\int_ {\bar \gamma} e^{-2\pi x/C}\,\dd x \,\dd \eta >0, 
    % \tag{\theequation a}
\end{equation}
where  $\bar \gamma$ denotes the region bounded by the chc $\gamma$. For $b(\rho)$ we let $h(x) = \cos 2\pi x - 2\pi \frac{S}{C}(\pi C \eta ^2) $, and note that  the level set of the chc is $\pi C \eta^2 = e^{2\pi x/C}E_{\sad}-\frac{C\rho}{2\pi}+\frac{g(x)}{4\pi}$. Then some algebra, that can be simplified using the relation below \eqref{eq:def_g}, gives
\begin{equation}
\label{eq:b_chc}
    b(\rho)=\frac{2\pi}{C}\int_{ \gamma}(f(x)+S\rho) e^{-2\pi x/C} \eta \,\dd x= \frac{2\pi}{C}\int_{\bar \gamma}(f(x)+S\rho) e^{-2\pi x/C}\,\dd x \, \dd \eta,
\end{equation}
 %\expand{\[\left (\tfrac{2\pi}{C}h(x)-h'(x)\right ) = \frac{2\pi}{C}(\cos 2\pi x - 2\pi \frac{S}{C} (e^{2\pi x/C}E_{\sad}-\frac{C\rho}{2\pi}+\frac{g(x)}{4\pi}))\] \[-(-2\pi\sin 2\pi x -2\pi \frac{S}{C}(2\pi/C e^{2\pi x/C}E_{\sad}+ \frac{g'(x)}{4\pi} ) ) = \]  \[\frac{2\pi}{C}(\cos 2\pi x +C\sin 2\pi x) +S\rho +S(-\frac{g(x)}{2C}+\frac{g'(x)}{4\pi}))\] Now use that we defined $g$ such that  $\frac{1}{2C}g(x)-\frac{1}{4\pi}g'(x)=\sin 2 \pi x$ from what we get desired result. Note that the fact that E term vanishe in the integral was trivial (it was a constant once mutiplied by exponential in the integral respect eta and was integrated over contractible curve (same start and finish)). }
where 
\begin{equation}
f(x)=\cos 2\pi x + (C-S)\sin 2 \pi x.
\label{eq:f}
\end{equation}
The chc tends to the concatenation of the two 
rhcs as $\rho \rightarrow \rho_0$ thus $C^{-1}b(\rho)/a(\rho)\rightarrow C^{-1}b/a$ from \eqref{eq:tilde_alpha_N}, i.e.~the chc curves tend to the approximate location of the necklace as $\rho \rightarrow \rho_0$, which is consistent with the unfolding of the $N$ point.
%\expand{All comes from pontriargin so all aproximate locations (to first order). } 
\begin{figure}[htbp]
    \centering
        {\includegraphics[width=0.75\textwidth]{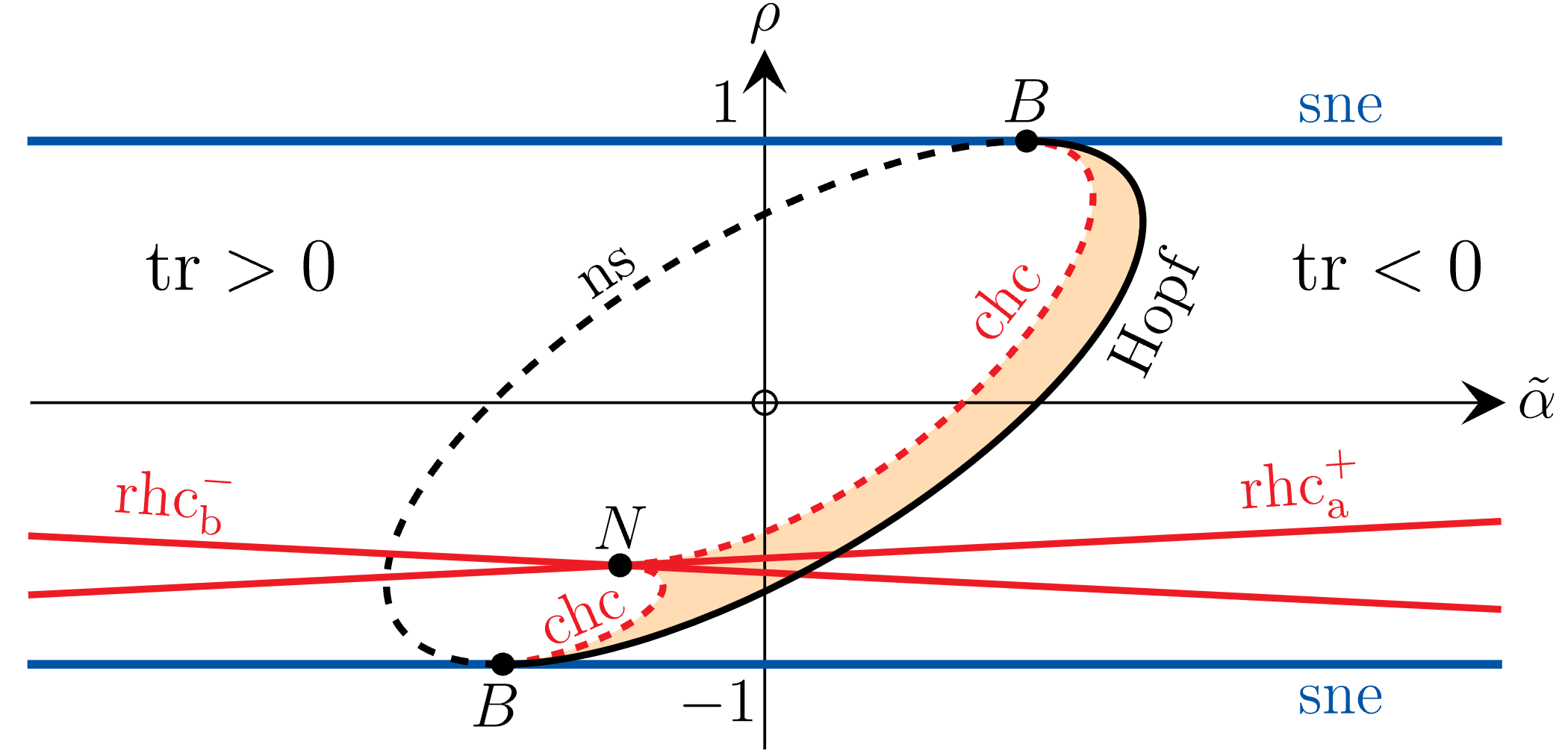}}
    \caption{Bifurcation diagram for \eqref{eq:reversible_approx_sqrt_eps} in $(\tilde \alpha, \rho)$ excluding the intersections of the rhc curves with ns and their consequences. The rhc curves are horizontal to first order but for clarity they and the chc curves are not drawn to scale. Parameter values in the orange shaded region have a cpo as will be shown in the following section.   }
    \label{fig:delta_neighbourhood}
\end{figure}

Note that $\tilde \alpha =C^{-1}b(\rho)/a(\rho)$, can be interpreted as the average value of  $(f+S\rho)/C$, in the region $\bar \gamma $, weighted by $e^{-2\pi x/C}$. When $\rho\rightarrow\pm 1$ the chc shrinks to a point so $\tilde \alpha$ tends to the value of $(f+S\rho)/C$ at the saddle. Thus, using $\sin 2\pi x_{\sad}=\rho$, we find $\tilde \alpha \rightarrow \pm 1 $ when $\rho\rightarrow\pm 1$.
%\expand{Since $\rho\rightarrow \pm 1$ and $\sin 2\pi x_{\sad}=\rho$ we must have cosinus tending to 0 and thus \[\frac{f(x_{\sad})+S\rho}{C}\rightarrow \frac{(C-S) (\pm1) +S (\pm 1)}{C} = \pm 1\]}
These are also the approximate locations of the $B$ points, which can be deduced from \eqref{eq:trace_zero_in_alpha_rho} and the fact that on the boundary of the resonance $\rho =\pm 1$. 
%\expand{We have $\rho ^2 + C^2 (\tilde \alpha -\rho) ^2 = 1 $ with $\rho =\pm 1$ then $C^2 (\tilde \alpha -\rho) ^2 =0$ so $\tilde \alpha = \rho = \pm 1$.     For approximate note that both the pontryagin method an the elipse is all to first order.}
Using the genericity of the $N$ and $B$ points we can conclude that the chc curves connect the $B$ points to the $N$ point as depicted in figure~\ref{fig:delta_neighbourhood}.
%\expand{First choose neighbourhoods $N$ and each $B$ point so that generic unfolding. Now can pick $\rho$ in compact intervals intersecting with this neighbourhoods, where we know that unique chc for each $\rho$. As this regions overlap they must represent same chc in overlaped region and thus we get connection. Only uniquely scenario is that in region of parameters missed by the intervals and neighbourhoods we have other chc curves, but doesn't change our conclusion directly.    }

Finally, we prove that the curves of chc do not intersect that of ns, so there is no $J$ point, justifying assumption 7a  and the genericity of the curves of chc.   Using $\sin 2\pi x_{\sad}=\rho$ and \eqref{eq:trace_zero_in_alpha_rho}, we find that the curve of ns has
\[\tilde \alpha \sim \frac{\cos 2\pi x_{\sad}+C\sin 2\pi x_{\sad} }{C} = \frac{f(x_{\sad}) + S\rho}{C}.\]
%\expand{From trace-zero equation, we get \[|\tilde \alpha -\rho| = \sqrt{1-\rho^2}/C\] so either $\tilde \alpha = \rho + \sqrt{1-\rho^2}/C $ or $\tilde \alpha  = \rho -\sqrt{1-\rho^2}/C$. The latter case is smaller so represents ns (other is Hopf which is to the right of ns). Now, as $\det = 4\pi^2 \cos 2\pi x_{\sad}$ is negative $\cos 2\pi x_{\sad } = -\sqrt{1-\rho^2} $ and we get the desired result.}
%NOTE: many/all times that we write sim we mean order sqrt or epsilon, (so metimes we go on to use this) 
Thus,  proving $\tilde \alpha_{\textrm{chc}}>\tilde \alpha_{\textrm{ns}}$ is equivalent  to showing that the average of $f$ over the region bounded by the chc (weighted by $e^{-2\pi x/C}$) is bigger than $f(x_{\sad})$, since the constant terms cancel out. 
%\expand{$\tilde \alpha_{\textrm{chc}}>\tilde \alpha_{\textrm{ns}}$ equivalent to average $\frac{f + S\rho}{C}$ bigger than value at saddle. We can cancel constant terms as  even in avarage they come out.  }
Showing this is tedious and not particularly insightful when $\phi\not =0$, but can be done; see appendix \ref{append:avarage_f}.

\section{Contractible periodic orbits}
\label{sec:cpo}
In this section we show numerically that the cpos of the reversible approximation give rise to a unique cpo for each parameter value in the zone between the curves of chc and the Hopf curve, i.e.~the shaded orange region in figure~\ref{fig:delta_neighbourhood}.
 %\expand{ Showing that chc to the left of Hopf  theoretically is probably hard as they are much closer (and tangent at B) than the ns curve with the chc. Note that showing that chc and ns didn't intersect took significant effort (no $J$ point). However, we do understand this region on $E$, $\rho$. In this parameters we know that the region is from from $E_{\min}$ to $E_{\max}$ (which is larger), and the first corresponds to Hopf and the later to chc, see argument below in the section. Then, the numeric monotonicity of $\tilde \alpha$ on $E$ allows us to get the result in $\tilde \alpha$, $\rho$.   } 
 Our analysis also shows that no cpo persist for parameters outside this region.  Moreover, the phase diagram of the reversible approximation does not allow for the appearance of any other cpo. In particular, there is at most one cpo for parameter values in a $\delta$-neighbourhood of the top or bottom of the resonance region. Outside these neighbourhoods we have shown in section \ref{sec:flow_|omega_y|<1-delta} and  \ref{sec:lack_of_chc_cpo} that there are no cpos, so we can conclude that assumption 7b is satisfied.
 %NOTE: can be interpreted as outside in general or outside but in R. Both are fine, just need, to do mental argument of only important part is R in different places.  
 %\expand{Note that we can only have cpo when we have equilibria so only for parameters in the resonance region which are all covered in section \ref{sec:flow_|omega_y|<1-delta} and  \ref{sec:Omega_y_mp_1|_leq_delta}. }
 %\expand{Didn't need to do this for chc as only need to control where they intersect ns (for $J$) points, and we know ns only in delta neigh.  }

We note that when applying Pontryagin’s method in section \ref{sec:chc} we never used that the specific level curve contained the saddle, nor its specific energy $E_{\sad}$.
%\expand{We need closed curve but that is general fact for Pontryagin as we explained. We did use contractibility to simplify computation (could be argued without maybe) and obviously to have integration of area bounded. However, it is clear enough that chc and cpo are contractible so we don't mention it. }
So exactly the same arguments show that a cpo $\gamma$ with energy $E$ persists at parameter value $\rho$ for a unique $\tilde \alpha $  given to first order by $\tilde \alpha (E, \rho)= C^{-1}b(E,\rho)/a(E,\rho)$, where from \eqref{eq:a_chc} and \eqref{eq:b_chc} we get 
\begin{equation*}
    \begin{split}
    a(E,\rho)&=\frac{4\pi}{C}\int_{x_{\min}}^{x_{\max}} e^{-2\pi x/C}\eta (x) \,\dd x,\\
    b(E,\rho)&=\frac{4\pi}{C}\int_{x_{\min}}^{x_{\max}}(f(x)+S\rho )e^{-2\pi x/C}\eta (x)\,\dd x,
    \end{split}
    \label{eq:a_b_cpo}
\end{equation*}
%\expand{Extra factor 2 comes from integrating the top and the bottom. As $\eta$ changes sign and direction of integration also flips, they compensate so we get twice the value on the top.  } 
with $x_{\min}$, $x_{\max}$ the intersections of $\gamma$ with the $x$-axis, and $\eta (x)$ the part of  $\gamma$ above it. If we denote by $E_{\cen}$ the energy at the centre inside the cpo and by $E_{\sad}$ the energy at the saddle on the chc enclosing the cpo, we have $E_{\cen}\leq E_{\sad}$ and  $E\in (E_{\cen}, E_{\sad})$.
%NOTE: there is a bit of ambiguity if energies are at center in unperturbed or perturbed. As everything is approximated, I think both interpretations would be fine, but we are thinking on the unperturbed case. 
%\expand{First note that max min energy in region bounded by chc has to be in boundary (chc) or in inflection point of $H$, that is equilibria system. To decide what is max and what min note that it is a general fact that in a Hamiltonian (this is after doing equivalence by exponential) centers are max/min energy function and saddles are saddles energy function (determinant vectorfield is same as determinant hessian). Then to determine if center is min or max note that it is at $\eta=0$ and when $|\eta|$ increases $H$ increases, so it must be minimum and maximum in chc. }
When $E\rightarrow E_{\sad}$, the cpo tends to the chc, so we get $\tilde \alpha (E,\rho)\rightarrow \tilde \alpha (\rho)$ from section  \ref{sec:chc}, the approximate value of $\tilde \alpha $ at the chc curves.
%\expand{We showed its genericity by pontryagin (so 1 dim manifold of them) and by no ns (which guarantees normal breaking). This was shown in no $J$ point, and was mentioned there. }
 Note that,  as in the case of chc, $\tilde \alpha (E,\rho)$ can be interpreted as the average of the function $(f+S\rho)/C$ (weighted by $e^{-2\pi x/C}$) in the region bounded by the cpo. So when  $E\rightarrow E_{\cen}$ and the cpo shrinks to a point,  $\tilde \alpha (E,\rho) $ tends to the value of the function at the centre. Using that $\sin 2\pi x_{\cen} =\rho$ we find that as $E\rightarrow E_{\cen}$,
\[\tilde \alpha (E,\rho) \rightarrow \rho +C^{-1}\sqrt{1-\rho^2},\]
%\expand{\[\frac{f(x_{\cen}) + S\rho}{C} = \frac{\cos 2\pi x_{\sad} + C\rho}{C} \] and then use that $\cos$ is positive in center as det is positive.}
which from \eqref{eq:trace_zero_in_alpha_rho}, is the approximate value of $\tilde \alpha$ at the Hopf curve. 
%\expand{From trace-zero equation, we get \[|\tilde \alpha -\rho| = \sqrt{1-\rho^2}/C\] so either $\tilde \alpha = \rho + \sqrt{1-\rho^2}/C $ or $\tilde \alpha  = \rho -\sqrt{1-\rho^2}/C$. The former case is bigger so represents Hopf (is is to the right of ns). }

Now if we show that for a fixed $\rho\in (-1,1)$, $\tilde \alpha (E,\rho)$ is strictly decreasing\footnote{Actually we need to show monotonicity for the exact value of $\tilde \alpha$. However, this can also be deduced from $\tilde \alpha'(E,\rho)$ as long as the error terms preserve their order upon differentiation. This follows from Pontryagin's method for parameters away from the boundary of the region of interest. Close to the boundary, uniqueness of cpo can be established from the unfolding of generic Hopf, chc, $N$ and $B$ points. } in the interval  $(E_{\cen}, E_{\sad})$, then we can conclude that for each parameter in the region  between the chc curve and the Hopf curve, we have exactly one cpo.
%NOTE: If we are using approximate energies (that is take the from unperturbed case) even if we compute this perfectly we need to use genericity of curves to cover region close to chc or Hopf.  

To show that $\tilde \alpha (E,\rho)$  is decreasing on $E$, we compute its derivative, 
%\expand{Actually method probably ok in Hopf. For close to chc, we know from pontryagin (no differentiation) that there can only be cpo close to conection. There generic unfolding chc implies only 1 created. }
\begin{equation}
    \label{eq:tilde_alpha_prime}
    \tilde \alpha'(E,\rho):=\frac{\partial \tilde \alpha}{\partial E} (E,\rho) = \frac{b'a -ba'}{Ca^2}
\end{equation}
where everything on the right is evaluated at $(E,\rho)$, and
\begin{equation*}
    \begin{split}
        a'(E,\rho)&:=\frac{\partial a}{\partial E}(E,\rho)=\frac{2}{C^2}\int_{x_{\min}}^{x_{\max}}\frac{1}{\eta(x)}\,\dd x,\\ 
        b'(E,\rho)&:=\frac{\partial b}{\partial E}(E,\rho)= \frac{2}{C^2}\int_{x_{\min}}^{x_{\max}}
        \frac{f(x)+S\rho}{\eta(x)}\,\dd x,
    \end{split}
\end{equation*}
which have been computed by differentiating under the integral sign, while using that $\eta(x_{\min})=\eta(x_{\max})=0$.  When computing  \eqref{eq:tilde_alpha_prime} one can remove the $S\rho$ term in $b$ and $b'$ using that $\tilde \alpha '=(\tilde \alpha-S\rho)'$.

As we are not able to compute the integrals, nor derive useful bounds on \eqref{eq:tilde_alpha_prime}, we check the sign of $\tilde \alpha'$ numerically.

\begin{figure}[htbp]
    \centering
        {\includegraphics[width=0.72\textwidth]{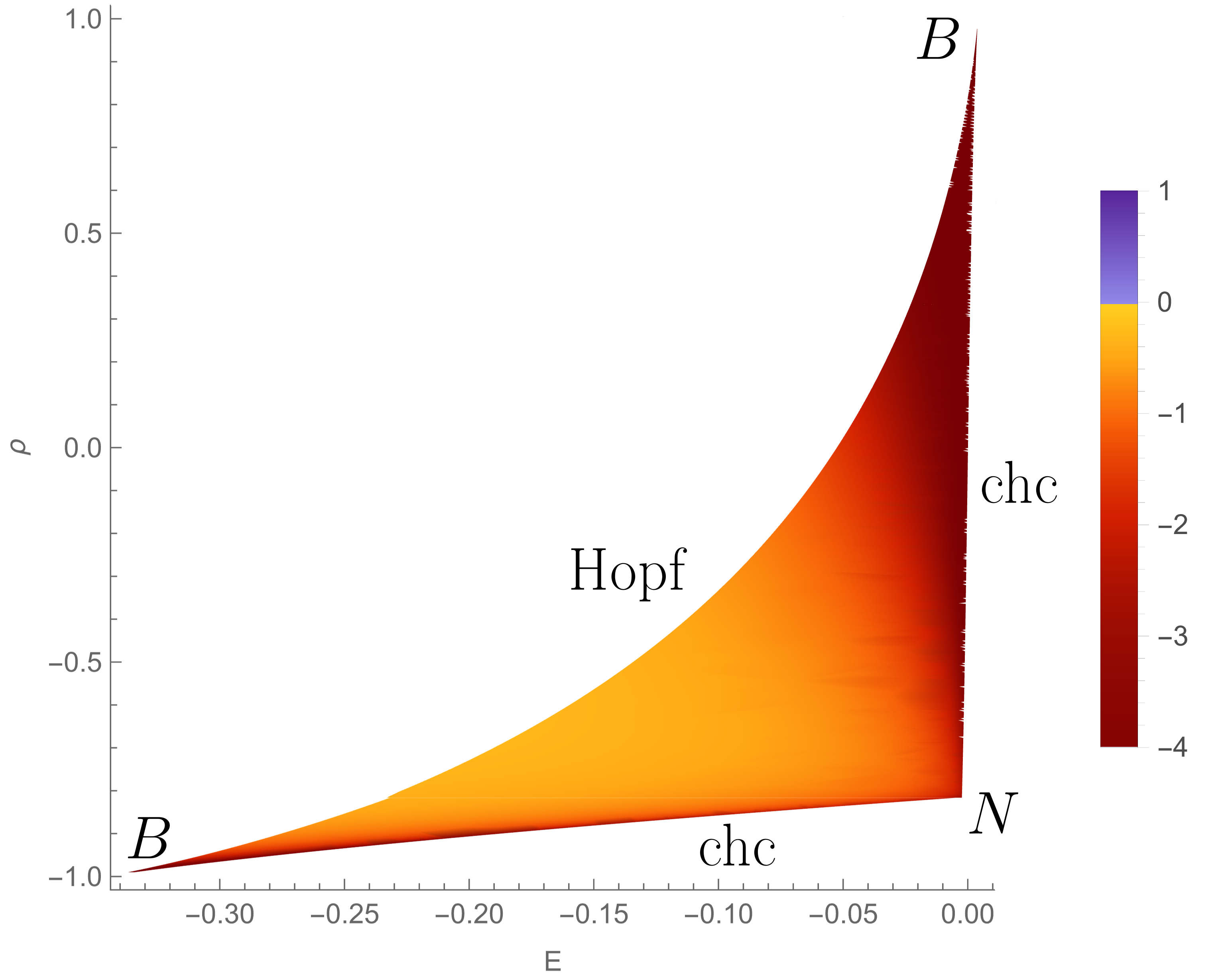}}
    \caption{ Density plot of $\tilde \alpha '$  in the regime of interest with $\phi=1/8$. There is no parameter value shaded with positive colour. Note that this region corresponds to the orange shaded region in figure~\ref{fig:delta_neighbourhood}. Note also that the righthand boundary in this figure is not quite vertical. }
    \label{fig:numerical_density_plot}
\end{figure}

To present the numerical computations we chose $\phi= \frac{1}{8}$, which is enough to make an example of the simplest bifurcation diagram. We note though, that we have not seen significant difference with other values $\phi \in ( \frac{1}{24}, \frac{5}{24})$. Figure~\ref{fig:numerical_density_plot} depicts a \texttt{Mathematica} density plot of $\tilde \alpha '$ in the region bounded by the Hopf curve and the chc curves. As the density plot does not have blue regions, or even light yellows, we are fairly confident that the value of $\tilde \alpha '$ is always negative and significantly away from zero. We can conclude then that $\tilde \alpha $ is decreasing on $E$.

To end this section we note that in the density plot it seems that $\tilde \alpha '$ has singularities at the $B$ points; this is expected for generic $B$ points. 

%NOTE: Numerically it seems that the minimum is reached at the center curve, with $\rho\approx -0.68$, where $\tilde \alpha '\approx -0.3$. 

\section{\texorpdfstring{$K$}{K} and \texorpdfstring{$H$}{H} points}
\label{sec:K_H_points}

Finally, we find the intersections of $\textrm{rhc}^\pm$ curves with the ns and the horizontal snp curves, which give rise to $K$ and $H$ points respectively. 

\begin{figure}[htbp]
    \centering
        {\includegraphics[width=0.99\textwidth]{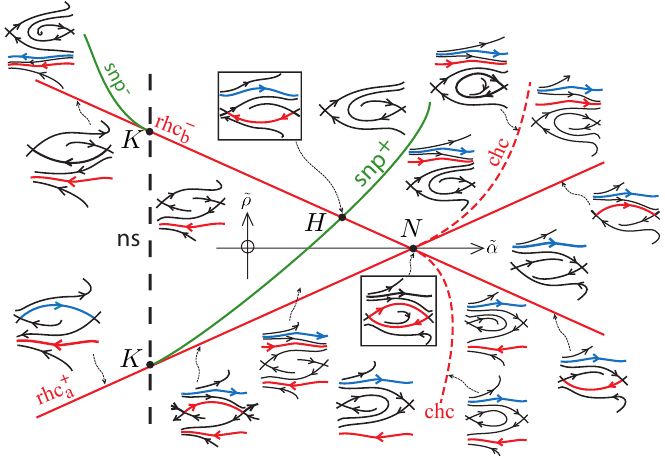}}
    \caption{   Bifurcation diagram in $(\tilde \alpha,\tilde \rho)$ coordinates around the $N$, $K$ and $H$ points leaving out the Hopf curve and the fan of rhc emanating from the $H$ point. The $\textnormal{snp}^{\pm}$  coincide to first order with the $\textnormal{rhc}^\pm$ curves but for clarity we have exaggerated their separation ($\pm$ indicates right-going/left-going).  }
    \label{fig:H_point}
    %NOTE: could depict the fan of rhc from H but I think better to leave it out. Should we check if rhc are drawn relatively similar to actual value with phi=1.1/24?
\end{figure}

Recall from \eqref{eq:rhc_rho} that the  $\textrm{rhc}^\pm$ curves are $\mathcal C ^1$-close to lines with steepness $\pm Ca/c$ in the $(\tilde \alpha, \tilde  \rho )$ plane.
%\expand{Get $C^1$ from pontriagin (we didn't mention it at the time though). Also not sure if true for all plan or only on compact set, but in any case that is enough. }
In this plane $\rho\sim \rho_0$, so that from \eqref{eq:trace_zero_in_alpha_rho}, we deduce that the ns curve has $\dd \tilde \alpha /\dd \rho \sim -S^2/C^2$.
%\expand{As ns is to the left of Hopf and using trace-zero equation, we have $\tilde \alpha = \rho - C^{-1}\sqrt{1-\rho^2}$. Differentiating (error terms behave well as we are $\mathcal C^1$ close to elipse) and evaluating at $\rho_0$ we get,      \[\dd \tilde \alpha /\dd \rho \sim 1- C^{-1}\frac{-2\rho_0 }{2\sqrt{1-\rho_0^2}} =1+ \frac{\rho_0 }{C\sqrt{1-\rho_0^2}} = 1+\frac{-(1+C^2)^{-1/2}}{C\sqrt{C^2/(1+C^2)}} = 1-C^{-2} = \frac{C^2-1}{C^2}=-\frac{S^2}{C^2} \]}
It must then be nearly vertical in the $(\tilde \alpha, \tilde  \rho )$ plane, due to the $\sqrt{\varepsilon}$ scaling, and thus it intersects once with each $\textrm{rhc}$ curve, see figure \ref{fig:H_point}.
%\expand{Note that cannot have any other intersection of this curves, as ns is contained in $\tilde \alpha$ and for this parameter we only have rhc in the region we are studying here. We have \[\frac{\dd \tilde \alpha }{\dd \tilde \rho}= \frac{\dd \rho }{\dd \tilde \rho}\frac{\dd \tilde \alpha }{\dd  \rho}\sim \sqrt{\varepsilon}(-S^2/C^2)\] then as error terms behave well we have transversal intersection so only one for each. }
This gives precisely two $K$ points, which are generic,  as the intersections are transverse and the trace has non-zero derivative across the ns curve.
%\expand{Essentially need transverse intersection of generic chc curve and generic ns curve. Generic chc is given by monotonicity slope in BM18. Generic ns is given by nonzero derivative. (Also need nothing else going on, ie. not intersection of 3 curves at same point, that is why for other cases we check trace non-zero for instance).  }
Their location can be found by substituting  $\rho\sim \rho_0$ in \eqref{eq:trace_zero_in_alpha_rho} and then using \eqref{eq:rhc_rho}, which gives 
\[\tilde \alpha \sim2\rho_0, \hspace{1cm}\tilde \rho \sim \mp\frac{16C(2-SC+4C^2)}{(1+C^2)^{3/4}(4+C^2)(4+9C^2) \tanh \frac{\pi}{C}}.\] 
%\expand{As before in the ns we have $\tilde \alpha \sim \rho -C ^{-1}\sqrt{1-\rho^2}$. Evaluating at $\rho_0$ \[\rho_0 -C^{-1}\frac{C}{\sqrt{1+C^2}}= 2\rho_0.\] Now plug this $\tilde \alpha$ in to rhc equation found by pontriagin to get $\tilde \rho$ (do it with computer). }
Putting together the top and the bottom of the resonance region proves assumption 8.

Generically, a $K$ point produces an arc of generic snp, tangent to the curve of rhc at the $K$ point. For this example, they must occur in the directions shown in figure \ref{fig:H_point} for compatibility with the sign of trace at the saddle, the direction in which the rhcs break and the rotational periodic orbits that exist below the necklace and not above, as illustrated in the figure.
The arcs of snp are $\mathcal C^1$-close to the corresponding $\textrm{rhc}$ curves in the $(\tilde \alpha,\tilde \rho)$ plane, as we will show in the upcoming section.  Thus, we have a unique intersection between snp and rhc curves of opposite direction, i.e.~an $H$ point, which is close to the $N$ point as depicted in figure \ref{fig:H_point}. Putting together the top and the bottom of the resonance region proves assumption 9. 
%\expand{Clearly snp- can intersect anything else in neighbourhood. As rhc+ has transverse intersection with rhc-, same is true fro snp+ (need $C^1$ close). In rest of delta neighbourhood we dont have more snp as reversible approximation doesn't allow such (note can only say stuff about horizontal as reversible approximation covers only some region of phase space, that is why we can have here non-horizontal rhc tha we do not detect on reversible approx). Outside delta neighbourhood it has already been studied and concluded that snp is on top of rhc (thus without intersections)}
%NOTE: outside delta neighbourhood already studied before and we said snp above rhc (so no intersection). Unclear if argument is 100\% rigorous there (because of possible existance of other snp curves). Regardless can't have this curves here intersecting outside, as they are above below rhc with same sign (going to other side delta neigh is impossible as on top of it homological direction is not horizontal, and inside delta neighbourghood, snp can only exist close to rhc as it is where reversible approximation allows such).

Moreover, the $H$ point is generic as its intersection is transverse and the trace at the saddle is non-zero on it. In particular, we have tongues emanating from the top of $H$ with rhc curves as boundaries and unique rotational periodic orbits (attractive as the trace is negative) for parameter values in their interior \cite{Baesens_2018}. In fact, there is a fan of tongues, one for each homology direction between $(\pm 1,0)$, where the homology direction of each tongue is given by the homotopy type and direction of its rotational periodic orbits. 
%\expand{We haven't define what homological direction is. We take it to be intuitive, homology of periodic orbit, take direction as a vector (signed direction actually). We don't use rational direction on purpose as irrational wouldn't make sense in this context (of periodic orbits)}
%NOTE: need direction to discard no primes eg: (2,2) but also to make sense of what in beteen means. 
This has been omitted from figure \ref{fig:H_point}, but it is depicted in figure \ref{fig:full_bif_diag}, where it is also shown how some of these tongues may extend outside a neighbourhood of $H$. 
%Note: Directions we get are positive (ie half plane is positive one, we go upwords) as outside neighbourhood dot y>0. 

\subsection{Transit map}
\label{sec:transit_map}
In this section we show that the snp curves are  $\mathcal C^1$-close to the corresponding $\textrm{rhc}$ curves in the $(\tilde \alpha,\tilde \rho)$ plane. In a neighbourhood of the $K$ points this follows from their genericity, see \cite{Perko_1994}. Outside this neighbourhood we use the transit map and without loss of generality we focus on the $\textrm{snp}^+$ case.

The saddle is near $x_g+\frac{1}{2}$ so it is natural to consider the transit map $\mathcal T$ from  $x=x_g$ to $x=x_g+1$. 
%\expand{The saddle is here as we are in $\tilde \rho$ so $\rho\sim \rho _0 $ and then then from the equation of the unperturbed rhc \eqref{eq:necklace_unperturbed} it is clear that we have the saddle when $\cos \pi (x-x_g) $ vanishes.  }
This map will be defined for $\eta>\eta_s$, the value of $\eta$ on the stable manifold at $x_g$, and will be mapped to $\hat \eta>\eta_u$, the value on the unstable one at $x_g+1$.
%\expand{Stable and unstable respect the affromentioned saddle at $x_g+1/2$. Note that this is the conection with positive $\eta$ which corresponds to rhc+. }
We will use the coordinate $z=\eta^2$ as $\eta_s, \eta_u>0$ and it will simplify some formulae.  

Using the Hamiltonian \eqref{eq:Ham} with $\rho = \rho_0$ we find that in the unperturbed case, where we have a necklace point, the transit map is simply
\begin{equation} 
\label{eq:unperturbed_mathcal_T}
\mathcal T(z)=e^{2\pi /C}(z-K),\end{equation}
with $K=(e^{-2\pi/C}-1)\rho_0/\pi^2>0$ and fixed point $z_s=z_u=-\rho_0/\pi^2>0$.
%\expand{We note that at $\rho_0$ we have $g(x) = -2\rho_0 C \cos 2\pi (x-x_g)$. Thus, \[H(x,\eta) = e^{-2\pi x /C}\left ( C \pi \eta ^2 +\frac{C\rho_0}{2\pi}(1+\cos 2\pi (x-x_g))\right ). \]  Now in the transit map we want same energy for $x_g$ and $x_g+1$, where the cosine is 1, so \[H(x_g,\eta) = H(x_g+1,\hat \eta ) \iff e^{-2\pi x_g /C}\left ( C \pi \eta ^2 +\frac{C\rho_0}{\pi}\right ) = e^{-2\pi (x_g+1) /C}\left ( C \pi \hat\eta ^2 +\frac{C\rho_0}{\pi}\right )\]   \[\hat \eta ^2 = e^{2\pi /C}(\eta ^2 +\frac{\rho_0}{\pi^2})-\frac{\rho_0}{\pi^2} =e^{2\pi /C}  \left (\eta ^2 -(e^{-2\pi /C}-1)\frac{\rho_0}{\pi^2}\right )\]  and replacing $z$ gives the desired map. Now as we are in the unperturbed case (with a necklace point) we know that the fixed point of the return map is given by the stable and unstable manifold, i.e we have rhc. To find the value of $z_u= z_s$ solve $\mathcal T (z) = z$.  }
%\expand{For the fixed point note that map only defined above $\eta_s$ and $\eta_u$ so in theory we don't include this values (the return map makes no sense for this values as the orbit converges to the saddle). However, if we take the above formule as a map it is clear that the affromention point is a fixed point and is just in the boundary of definition. More generaly, the transit map can always be exteded continuously by $T(z_s) = z_u$ as orbits starting close to $\eta_s$ will end up close to $\eta_u$  (in general no differentiable as slope can be infinity).  }
In the perturbed case, $z_u=z_s$ iff there exists a $\textnormal{rhc}^+$. So the rhc is detected by
\begin{equation}
    \label{eq:z_u-z_s=0}
    z_u-z_s=0
\end{equation}
where by the stable manifold theorem $z_u$ and $z_s$ are smooth on parameters.
%\expand{Maybe also need to use smooth dependence respect parameters of solutions ODEs. More concretely we use stable manifold theorem to get away from neigh equilibrium (where time goes to infinity). Then smooth dependence allows us to get to the desired transverse section.  }

From \eqref{eq:unperturbed_mathcal_T} we note that the slope of the perturbed transit map is close to $e^{2\pi /C}>1$ for $z$ sufficiently above $z_s$, because the return time is bounded, so there is no snp there. Thus, it is left to study  $\mathcal T'(z)$ when $z$ is close to $z_s$, which we will quantify shortly, where the return time goes to infinity. 

We will study $\mathcal T'(z)$ by exploring the evolution of the symplectic area $\mathcal A=e^{-2\pi x/C}  \dd x\wedge\dd\eta $. First consider the coordinate  $\chi = -\frac{C}{2\pi} e^{-2\pi x/ C }$ so that $\mathcal A$ is the standard area, i.e. $\mathcal A =   \dd \chi\wedge\dd\eta  $. Then, its rate of growth is given by the divergence of the vector field in these coordinates, that is $\overset{\hspace{2.8pt}\sbullet}{\mathcal  A}= \tau\mathcal A $, where  
\begin{equation}
\label{eq:dot_A}
\tau  = 2\pi \sqrt{\varepsilon}\left (\sin 2\pi x -\tilde \alpha +\frac{1}{C} \cos 2\pi x -2\pi ^2\frac{S}{C}\eta^2\right ) +O(\varepsilon).\end{equation}
%\expand{To compute this as described above, note that $\primebullet {\mathcal A} = \textnormal{div}(\primebullet \chi, \primebullet \eta) \mathcal A$. To compute divergence note that $\primebullet \chi = e^{-2\pi x(\chi) /C}\primebullet x $ and then compute derivatives only using that $x = \frac{-C}{2\pi}\log(\frac{-2\pi}{C} \chi)$ when necessary.   }
%\expand{An alternative is to change time (advantage that then clear first order cancels out). First let's denote by $\tilde {\mathbf{v}}$ the equivalent vector field $e^{-2\pi x /C} (\primebullet x,\primebullet \eta)$ and recall that its first order approximation is Hamiltonian. Let  $\mathrm t$ be its associate time so that     \[\frac{\dd}{\dd \mathrm t} x =\frac{\dd s }{\dd \mathrm t} \frac{\dd}{\dd  s} x = e^{-2\pi x /C} \primebullet x  = \tilde v_1.\]    Then if we denote the standard area by $\omega = \dd x\wedge\dd\eta$ we have, \[\primebullet {\mathcal A } = \frac{\dd}{\dd \mathrm t}\left (\frac{\dd \mathrm t}{\dd s} e^{-2\pi x /C} \omega  \right ) = \frac{\dd}{\dd \mathrm t} \omega = \textnormal{div} (\tilde v) \omega = ( e^{2\pi x /C}\textnormal{div} (\tilde v) ) \mathcal A,\]  where we have used that the growth rate of the standard area, $\omega$, is given by the divergence of the vector field.  }
So if we consider the path $\gamma$ for the flow from $(x_g,\eta)$ to $(x_g+1,\hat \eta )$ we have 
\begin{equation*}
%\label{eq:A/hat A}
\hat{ \mathcal A}= \mathcal A e^{\int_\gamma \tau \, \dd s}.
\end{equation*}
Applying $\mathcal A$ to the vector field \eqref{eq:reversible_approx_sqrt_eps} and a vertical vector $\delta \eta$ at the initial point yields $e^{-2\pi x_g/C} (2\pi C \eta + O(\sqrt{\varepsilon}))\, \delta \eta$.  
During the transit, this area gets multiplied by $e^{\int_\gamma \tau\, ds}$.  The vector field is carried to the vector field at the image point; the vertical vector is carried to a vector that in general is not vertical, but after adding a suitable multiple of the vector field it can be made a vertical vector $\delta\hat{\eta}$, spanning the same area. This pair spans area $e^{-2\pi (x_g+1)/C} (2\pi C \hat{\eta} + O(\sqrt{\varepsilon}))\, \delta \hat{\eta}$.
 So this yields
\[\mathcal T'(z)=\frac{\delta (\hat \eta ^2)}{\delta ( \eta ^2)}= e^{\int_\gamma \tau \,\dd s} \left (e^{2\pi/C}+O(\sqrt{\varepsilon})\right ),\]
as $\eta > \eta_s$, $\hat \eta > \eta_s$ and $\eta_s, \eta_u\sim \sqrt{-\rho_0}/\pi >0$ % $\eta_s, \eta_u\sim  \pi ^{-1}(1+C^2)^{-1/4}>0$ 
which has been deduced from \eqref{eq:unperturbed_mathcal_T}.
%\expand{\[\mathcal T'(z)=\frac{\delta (\hat \eta ^2)}{\delta ( \eta ^2)}= \frac{2\hat \eta \delta \hat \eta }{2\eta \delta \eta } \] and \[e^{\int_\gamma \tau \,\dd s}  = \frac{\hat{\mathcal A}}{\mathcal A} = \frac{e^{-2\pi (x_g+1)/C} (2\pi C \hat{\eta} + O(\sqrt{\varepsilon}))\, \delta \hat{\eta}}{e^{-2\pi x_g/C} (2\pi C \eta + O(\sqrt{\varepsilon}))\, \delta \eta} =\frac{2\hat \eta \delta \hat \eta }{2\eta \delta \eta } (e^{-2\pi/C}+O(\sqrt{\varepsilon}))\] where in the last equality to deal with the error terms we have used that $\eta$ and $\hat \eta $ are not close to zero as they are bounded by a constant. }

To compute the integral in the exponent consider a neighbourhood $U$ of the saddle where a smooth conjugacy flattens the stable and unstable manifolds to coincide with the axes. We can choose $U$ independent of $\varepsilon$ and such that $\gamma$ intersects the boundary transversely for $\eta$ close to $\eta_s$.
%\expand{Here by close we mean close with a constant independent of $\varepsilon$ as neigh also chose independent of $\varepsilon$. Thus, it is clear that any close in terms of $\varepsilon$ will be much smaller (if we choose $\varepsilon$ small enough). So that the regin of $\eta$ that we are studying intersects with this.   }
%NOTE: also should choose it such that it moves with parameters smoothly. This may not be needed as parameter affect only when eps\not =0, but still we net the conjugation to change smoothly so constants below move smoothly
Then, the time spent outside $U$ is bounded so the contribution of the integral there is $O(\sqrt{\varepsilon})$.
%\expand{Bound integral by time (bounded) times max trace that as $\eta $ close to $\eta_s$ term in parenthesis trace is bounded and so is order $\sqrt{\varepsilon}$. }
The contribution in $U$ can be computed after flattening of the invariant manifolds as is done in appendix~\ref{append:integral_trace_saddle}. 
 In total we get, 
%NOTE: cant simply say value at saddle times time, as neighbourhood is not small enogh.  
\begin{equation}
    \int_{\gamma} \tau \,\dd s= T \tau_{\sad}+O(\sqrt{\varepsilon})=T (\sqrt{\varepsilon}\,\tau_0+O(\varepsilon))+O(\sqrt{\varepsilon}),
    \label{eq:int_tau}
\end{equation}
where $T$ is the time spent in $U$ and in the second equality we approximate to first order $\tau$  at the saddle, $\tau_{\sad}$, by evaluating  \eqref{eq:dot_A} at $(x_g+\frac{1}{2},0)$ so, 
\[\tau_0=-2\pi \,(2\sin 2\pi x_g+\tilde \alpha).\]
If we let $\lambda$ be the unstable eigenvalue, $\lambda \sim  2\pi C(1+C^2)^{-1/4}$,
%\expand{In the unperturbed trace is zero so eignevalue is given by $\sqrt{-\det}$. Now $\det = 4\pi^2C\cos 2\pi x$ and at equilibria $\sin 2\pi x = \rho_0$. Then substituding $\cos 2\pi x = -\sqrt{1-\rho^0}$ we get the result (sign as in saddle trace is negative). }
the time $T$  spent near the saddle is given by, $-\lambda^{-1} \log m (\eta -\eta_s)$ where 
$m$ is around 1, that is, $m=O(1)$ and $1=O(m)$. Putting everything together, 

\begin{equation} 
\mathcal T'(z)=\big (e^{2\pi/C}+O(\sqrt{\varepsilon}) \big ) \left (m(\eta -\eta_s )\right )^{\sqrt{\varepsilon}\,\nu+O(\varepsilon)} , 
\label{eq:T'(z)_eta}
\end{equation}
where $\nu= -\lambda ^{-1}\tau_0>0$ and is around 1 as we are to the right of the neutral saddle curve (where the trace is zero) and outside a neighbourhood of the $K$ point. 
%\expand{Note that around one thinking that in the $\tilde \alpha =O(1)$ otherwise, we are essentially leaving the $\delta$ neighbourhood.}
 We are interested in points with slope 1, so that the second factor in \eqref{eq:T'(z)_eta} must be significantly smaller than one, in particular $\eta = \eta_s + O(\sqrt{\varepsilon})$.
%\expand{Note $\lim_{x->0} x^x=1$ (by taking log for instance and polinomial beat log) then if difference of etas is $\sqrt{\epsilon}$ or bigger, slope will be close to $e^{2\pi C}$ or bigger }
Thus, 
\[z-z_s  = (2\eta_s + O(\sqrt \varepsilon) )(\eta -\eta_s)\]
%\expand{\[z-z_s = (\eta + \eta_s)(\eta -\eta_s) = (2\eta_s + O(\sqrt \varepsilon) )(\eta -\eta_s)\]}
and converting to $z$ coordinates,
\begin{equation}
 \label{eq:mathcal_T'}   
\mathcal T'(z)=\left (e^{2\pi/C}+O(\sqrt{\varepsilon})\right )\left (\frac{z-z_s}{k+O(\sqrt{\varepsilon})}\right )^{\sqrt{\varepsilon}\, \nu+O(\varepsilon)} , \end{equation}
where $k=2\eta_s /m$  and is  around 1.
Then, points with slope 1 are given by,
\begin{equation}
    \label{eq:slope_1}
 z-z_s=\left (k+O(\sqrt{\varepsilon})\right )\left (e^{-2\pi/C}+O(\sqrt{\varepsilon})\right )^{\frac{1}{\sqrt{\varepsilon}\hspace{1pt}\nu+O(\varepsilon)}}.
 \end{equation}
 
We can compute the transit map by integrating $\mathcal T'$ while using that $T(z_s)=z_u$, so 
\begin{equation}
\label{eq:mathcal_T}
\mathcal T(z) = z_u + \left (k e^{2\pi/C}+O(\sqrt{\varepsilon})\right ) \left (\frac{z-z_s}{k+O(\sqrt{\varepsilon})}\right )^{1+ O(\sqrt{\varepsilon})}.\end{equation}
We note that some error terms of $\mathcal T'$ depend on $z$, so when we integrate treating them as constants we lose some control. In particular, one can check that all the errors in this section before \eqref{eq:mathcal_T} have derivatives with respect to parameters of the same order, whereas we cannot guarantee this for the terms in \eqref{eq:mathcal_T} which means we will have to be careful in our analysis.

Recall now that snp correspond to fixed points of $\mathcal T$ with slope 1, which by doing some algebra with \eqref{eq:slope_1} and \eqref{eq:mathcal_T}   occur only  when, 
\begin{equation}
\label{eq:z_u-z_s=exponential}
    z_u-z_s=O(1)\left ( e^{-2\pi /C}+O(\sqrt{\varepsilon}) \right )^{\frac{1+O(\sqrt{\varepsilon})}{\sqrt{\varepsilon}\hspace{1pt}\nu+O(\varepsilon)}}.
\end{equation}
%\expand{This equation give parameters for which we can (and indeed do have for some $z$) a fixed point with slope 1, i.e. a snp.  }
The right hand side of the expression above is exponentially small in $\sqrt{\varepsilon}$, so to all orders this equation is equivalent to \eqref{eq:z_u-z_s=0}, i.e.~the $\textnormal{snp}^+$ is indistinguishable from the $\textnormal{rhc}^+$.

However, we cannot conclude that the $\textnormal{snp}^+$ curve is $\mathcal C^1$-close to the $\textnormal{rhc}^+$ by simply differentiating \eqref{eq:z_u-z_s=exponential} with respect to parameters, as we do not have control over the derivatives of the error terms. So instead we consider the point $z$ with slope 1, defined by \eqref{eq:slope_1}. 
Then we can compute the derivative of $\mathcal T(z)$ with respect to parameters by differentiating under the integral sign. For instance with respect to $\tilde \alpha$ we have,
\[\frac{\dd}{\dd \tilde \alpha } \mathcal T(z) = \frac{\dd }{\dd \tilde \alpha }z_u + \frac{\dd}{\dd \tilde \alpha } z- \frac{\dd}{\dd \tilde \alpha }z_s +\int_{z_s}^{z}\frac{\dd}{\dd \tilde \alpha } \mathcal T'(\tilde z) \,\dd \tilde z ,   \]
so that differentiating $\mathcal T(z)=z$ we get,
\[\frac{\dd}{\dd \tilde \alpha}(z_u-z_s) = \int_0^{z-z_s}\frac{\dd}{\dd \tilde \alpha } \mathcal T'(\Delta z+z_s) \,\dd \Delta z .   \]
If we show that the right hand side is  exponentially small in $\sqrt{\varepsilon}$, this expression coincides to all orders with the derivative of \eqref{eq:z_u-z_s=0} and thus, the $\textnormal{snp}^+$ curve is $C^1$-close to the $\textnormal{rhc}^+$ curve in the region parametrised by $\tilde \alpha$, $\tilde \rho$. Bounding the right hand side is a bit convoluted. One first  differentiates \eqref{eq:mathcal_T'} with respect to $\tilde \alpha$ while using that the error terms preserve their order and that $\frac{\dd}{\dd \tilde \alpha }k,\frac{\dd}{\dd \tilde \alpha }\nu = O(1)$. The resulting expression can be integrated explicitly and then evaluated using \eqref{eq:slope_1} which leads to
\[\int_0^{z-z_s}\frac{\dd}{\dd\tilde \alpha }\mathcal T'(\Delta z+z_s)\,\dd\Delta z = O(1/\sqrt{\varepsilon}) \left ( e^{-2\pi /C}+O(\sqrt{\varepsilon}) \right )^{\frac{1+O(\sqrt{\varepsilon})}{\sqrt{\varepsilon}\hspace{1pt}\nu+O(\varepsilon)}}\]
%EXPAND: For differentiation use formula of $(f^g)'=f^g logf g' + f^g f'/f$ 
which is indeed exponentially small.

\section*{Acknowledgments}
We would like to thank David Marín for helpful
discussions both in person and through correspondence. We are also grateful for comments made by Armengol Gasull, Joan Torregrosa and Jordi Villadelprat in the  \textit{New Trends on Bifurcations in Ordinary Differential Equations} conference. Finally we would like to give our sincere thanks to the referee for an extremely thorough reading of our work and their many useful suggestions.  
%other people at uab from conference? 

\appendix

\section{Generic saddle-node curves}
\label{append:generic_saddle-node}

The saddle-node curves are given by $\det =0$,  but for this section we exclude the points with zero trace, because they have double eigenvalue zero so are not generic saddle-nodes.  At the rest, non-zero trace implies simple eigenvalue 0, which is one condition for genericity, and we show the remaining conditions.

Following  \cite{Kuznetsov_2004},
%REF(NONLINEAR OSSILATIONS THM3.4.1)
we can first verify the transversality by checking that the derivative of the map $(\vect x,\vect\Omega)\mapsto (\vect G, \det)$ (where $\vect G$ denotes the vector field) has full rank. As $D_{\vect \Omega} \vect G$ is the identity and $D_{\vect \Omega} \det=\vect 0$, it is enough to check that $D_{\vect x}\det\not =\vect 0$ at the saddle-nodes. 
%\expand{First note that     \[D \textnormal{map} = \begin{pmatrix}    D_x G_1 & D_y G_2 & 1 &0\\     D_x G_2 & D_y G_2 & 0 &1\\    D_x \det &D_y \det & 0 &0 \end{pmatrix} \] (where to check $D_{\vect \Omega} \det=\vect 0$ we are doing partial derivatives, so only checking if in expression of $\det$ the parameters appear. Don't consider secundary dependencies). So for maximal rank only need to check that $D_{\vect x}\det\not =\vect 0$, as then we will have 3 independent rows. We need to check this at the  saddle-node (there is where transversality condition needs to be verified), so at the $\det = 0$ curve.  }
Note that   $\frac{\partial \det}{\partial x}=0$ and $\det =0 $ implies that $\psi'(y)$ is orthogonal to an orthonormal basis, which is impossible as $\psi'(y)$ does not vanish, so the rank condition is satisfied.
%\expand{We need $\det = 0$ because we are at saddle-node. Now       \[\det =- 2\pi \varepsilon(\cos 2\pi x, \sin 2\pi x)\cdot \psi'(y),\]   \[\frac{\partial}{\partial x}\det =- 4\pi^2 \varepsilon(-\sin 2\pi x, \cos 2\pi x)\cdot \psi'(y)\]      And note that the vectors on $x$ are orthogonal. $\psi '(y)$ never vanishes as it gives another elipse.  }

Let $A = D_{\vect x} \vect G_{(\vect x,\vect \Omega^e(\vect x))}$ be the derivative of the vector field with respect to the first coordinate at $(\vect x,\vect \Omega^e(\vect x))$ and similarly $B= D^2_{\vect x} \vect G_{(\vect x,\vect \Omega^e(\vect x))}$.
Then, the non-degeneracy of the required second-order Taylor coefficient can be checked following  \cite[Theorem~3.4.1]{Nonlinear_oscillations_1983}, by
\[\vect p B(\vect q,\vect q)\not =0,\]
where $\vect q $ and  $\vect p$ are respectively a right and a left eigenvector of eigenvalue zero of $A$.

In our case $A$ is given by \eqref{eq:linearization_equilibria}, and using that at a saddle-node $\det=0$ we find,
\[\vect p=\begin{pmatrix}
\cos 2\pi x & \sin 2\pi x
\end{pmatrix}, \hspace{10pt}
\vect q=\begin{pmatrix}
d & - \varepsilon
\end{pmatrix}^T,\]
with $d=\frac{\cos 2\pi y}{\cos 2 \pi x}=\frac{\sin 2\pi (y-\phi)}{\sin 2\pi x}$, where the equality follows from $\det=0$ and  is always well defined as the denominators do not vanish simultaneously. 
We have 
\[B(\tilde {\vect x},\hat{\vect x})= 4\pi^2\begin{pmatrix}
\varepsilon \cos 2\pi x\hspace{4pt}\tilde x\hat x+\cos 2\pi (y-\phi) \hspace{4pt} \tilde y \hat y \\
\varepsilon \sin 2\pi x\hspace{4pt}\tilde x\hat x+\sin 2\pi y \hspace{4pt} \tilde y \hat y  \hphantom{\cos 2\pi } %to make + coincide in matrix
\end{pmatrix},\]
so 
\begin{equation}
\label{eq:coef_sn}
\vect p B(\vect q,\vect q)=4\pi^2\varepsilon (d^2 + \varepsilon(\cos 2\pi x\cos 2\pi (y-\phi)+\sin 2\pi x\sin 2\pi y)).\end{equation}
Note that $d$ can not be close to zero, as $\cos 2\pi y$ and $\sin 2\pi (y-\phi)$ do not vanish simultaneously. Thus, \eqref{eq:coef_sn} does not vanish and the saddle-node curves are generic.

\section{Non-degeneracy of \texorpdfstring{$B$}{B} points}
\label{append:generic_B}
First, we need to find the $B$ points to first order. As they are on the trace-zero loops, they have $y\sim \pm \frac{1}{4}$ and then substituting $y$ in $\det =0$ we find $x\sim \pm \frac{1}{4}$ (not necessarily with the same signs). 
%\expand{\[\det =- 4\pi^2 \varepsilon(\cos 2\pi x, \sin 2\pi x)\cdot (-C\sin 2\pi y  + S \cos 2\pi y, \cos 2\pi y)).\]     So evaluating at $y\sim \pm \frac{1}{4}$ we have    \[0\sim (\cos 2\pi x) (\mp C) \]      So to have this need $x\sim \pm\frac{1}{4}$ regarless sign in previous equation. If we evaluate $y=\pm\frac{1}{4}+O(\varepsilon)$ (as is what we actually have in \tr = 0, see section in B points in main text) we get  \[0=(\cos 2\pi x) (\mp C)+O(\varepsilon)\]   as  $\cos$ in $\frac{1}{4}$ has taylor of $\sin$ at 0 and thus x=\pm \frac{1}{4}+O(\varepsilon). }

Without loss of generality we consider the $B$ point on the outer saddle-node curve and near the top of the resonance region, i.e.~$x,y\sim \frac{1}{4}$. Similarly to the previous section, checking that the rank of the derivative of the map $(\vect x,\vect \Omega)\mapsto (\vect G,\det, \tr)$ is maximal, see  \cite{Kuznetsov_2004}, is equivalent to the independence of the vectors 
%\expand{First note that     \[D \textnormal{map} = \begin{pmatrix}    D_x G_1 & D_y G_2 & 1 &0\\     D_x G_2 & D_y G_2 & 0 &1\\    D_x \det &D_y \det & 0 &0 \\   D_x\tr &D_y \tr & 0 &0  \end{pmatrix} \] (where partial derivatives so no consider secondary dependencies). So for maximal rank check that rows independent. Note that $B$ points are actually at $x,y=\pm\frac{1}{4}+O(\varepsilon)$ as discussed in previous expand, and thus we get     \[\frac{\partial\det}{\partial x} =- 8\pi^3 \varepsilon(-\sin 2\pi x, \cos 2\pi x)\cdot (-C\sin 2\pi y  + S \cos 2\pi y, \cos 2\pi y) \sim - 8\pi^3 \varepsilon (-1,0) \cdot (-C,0) \]    \[\frac{\partial\det}{\partial y} =- 8\pi^3 \varepsilon(\cos 2\pi x, \sin 2\pi x)\cdot (-C\cos 2\pi y  - S \sin 2\pi y, -\sin 2\pi y) \sim - 8\pi^3 \varepsilon (0,1) \cdot (-S,-1)\]       \[\frac{\partial \tr}{\partial x} = 4\pi^2 (\varepsilon\cos 2\pi x) = O(\varepsilon^2)\]    \[\frac{\partial \tr}{\partial y} = 4\pi^2 (\sin 2\pi y) =4\pi ^2+ O(\varepsilon^2)\]     So we get desired results. Note however, that only knowing that $D_{\vect x}\tr  = 4\pi ^2(o(1),1+O(1))$ would be enough for our results and would could be computing only knowing that $x,y\sim \pm 1/4$.    }
\[D_{\vect x}\det \sim 8 \pi^3\varepsilon(-C,1),\hspace{1cm}D_{\vect x}\tr  = 4\pi ^2(O(\varepsilon^2),1+O(\varepsilon^2)),\]
%\expand{easier to get o(1) error term in tr and enogh for conclusion, but we give best bounds here. }
which have been evaluated at the $B$ point. They are clearly independent for $\varepsilon$ small enough.

For the non-degeneracy condition we let $s=\sin 2\pi x$ and keep the notation introduced in  the previous section of the appendix, so that 
\[A=2\pi s\begin{pmatrix}
\varepsilon & d\\
-\varepsilon^2 d^{-1} & - \varepsilon
\end{pmatrix}, \hspace{10pt} 
B(\tilde {\vect x},\hat {\vect x})=4\pi^2\begin{pmatrix}
\varepsilon \cos 2\pi x\hspace{4pt} \tilde x \hat x + \cos 2\pi (y-\phi)\hspace{4pt}\tilde y \hat y\\
\hphantom{o} \varepsilon s\hspace{4pt} \tilde x \hat x + \sin 2\pi y \hspace{4pt} \tilde y \hat y \hphantom{l} %\hphantom to make + coincide in matrix
\end{pmatrix},\]
as we are now  evaluating  at the $B$ points where both the trace and the determinant vanish.  
Let,
\[\vect q_0=\begin{pmatrix} d\\-\varepsilon \end{pmatrix},\hspace{20pt}
\vect q_1=\frac{d}{2\pi s( d^2+\varepsilon^2)}\begin{pmatrix} \varepsilon\\d \end{pmatrix}\]
being respectively a null vector of  $A$ ($A\vect q_0=0)$ and an independent vector scaled to make $A\vect q_1 =\vect  q_0$.  Let
\[
\vect p_1=2\pi s\begin{pmatrix} \varepsilon d^{-1} & 1 \end{pmatrix},\hspace{20pt}
\vect p_0=\frac{1}{d^2+\varepsilon^2}\begin{pmatrix} d& -\varepsilon \end{pmatrix},\]
being respectively a null form of  $A$ ($\vect p_1 A =0)$ scaled to make $\vect p_1\vect q_1 = 1$, and a form $\vect p_0$ chosen so that $\vect p_0\vect q_1 = 0$ and scaled to make $\vect p_0A = \vect p_1$ (equivalently $\vect p_0\vect q_0 = 1$).

Now, from \cite{Kuznetsov_2004}, the non-degeneracy condition is that the following second-order Taylor coefficients do not vanish 
\begin{align*}
    a&= \frac{1}{2}\vect p_1 B(\vect q_0, \vect q_0)=4\pi^3 \varepsilon C^2+O(\varepsilon^2),\\
    b&=\vect p_0 B(\vect q_0,\vect q_0)+\vect p_1 B(\vect q_0,\vect q_1)=-4\pi ^2\varepsilon+O(\varepsilon^2),
\end{align*}
%EXPAND: weird choice of where to put epsilon in a, but for consistency in general it was always to the left (so can see easely) except for the universal constants (usually 2 pi). 
where we have used that $s\sim 1$ and $d\sim C$ in the $B$ point with $x,y\sim \frac{1}{4}$. 
%EXPAND:  technically for computation  need 1+O... and also need need x=1/4+.. to evaluate cosinus (maybe)
Thus, for $\varepsilon$ small enough, the coefficients are non-zero and we can conclude that the $B$ point is generic. Note that a factor was missed in \cite{Baesens_2018}, which is corrected by the above results.
%\expand{I think essentially they missed the $2\pi$ in front of $A$. This causes a missing factor of $1/2\pi$ in $\vect q_1$ and a $2\pi $ factor missing in $\vect p_1 $ and $a$. }

\section{Lyapunov coefficient of the centres}
\label{append:lyapunov_coef}

As for the other bifurcations we check the transversality condition by noting that in the top/bottom of the resonance region, 
\[D_{\vect x}\tr = 4\pi^2(O(\varepsilon),\pm 1+O(\varepsilon ^2))\not =\vect 0,\]
so the differential of the map  $(\vect x,\vect\Omega)\mapsto (\vect G, \tr)$ has maximal rank. 
%\expand{First note that     \[D \textnormal{map} = \begin{pmatrix}    D_x G_1 & D_y G_2 & 1 &0\\     D_x G_2 & D_y G_2 & 0 &1\\    D_x\tr &D_y \tr & 0 &0  \end{pmatrix} \] (where partial derivatives so no consider secondary dependencies). So for maximal rank check that rows independent. Note that trace-zero is at $y=\pm\frac{1}{4}+O(\varepsilon)$ as discussed in previous $B$ points section main text, and thus we get     \[\frac{\partial \tr}{\partial x} = 4\pi^2 (\varepsilon\cos 2\pi x) = O(\varepsilon)\]    \[\frac{\partial \tr}{\partial y} = 4\pi^2 (\sin 2\pi y) =\pm 4\pi ^2+ O(\varepsilon^2)\]     So we get desired results. Note however, that only knowing that $D_{\vect x}\tr  = 4\pi ^2(o(1),1+O(1))$ would be enough for our results and would could be computing only knowing that $x,y\sim \pm 1/4$.    }

We now want to find the Lyapunov coefficient of the centres on the trace-zero curves to check that we have  non-degenerate Hopf bifurcations. 

We follow \cite{Kuznetsov_2004} by first finding an affine coordinate change, 
\begin{equation}
    \tilde{ \vect x}=\begin{pmatrix} a & b\\ c&d\end{pmatrix}\delta \vect x, 
    \label{eq:Lyapunov_coef:coord_change}
\end{equation}
where $\delta \vect x = \vect x-\vect x_0$, with $\Delta = ad-bc\not = 0$ and inverse
\[\delta \vect x=\frac{1}{\Delta}\begin{pmatrix} d & -b\\ -c&a\end{pmatrix}\tilde{ \vect x}, \]
to reduce the linearised dynamics around the centre $\vect x _ 0$ to the form, 
\begin{equation}
\label{eq:Lyapunov_coef:linear_normal_form}
\dot {\tilde {\vect x}}= \begin{pmatrix}
0 & - \omega\\ \omega & 0 
\end{pmatrix}\tilde {\vect x}.
\end{equation}
The equilibria with $\tr=0$ can be parametrised by $x_0\in \mathbb R$ and the choice of root $y_0$  of $\cos 2\pi y_0=\varepsilon \sin 2 \pi x_0$ close to $ \frac{1}{4}$ (resp. $ -\frac{1}{4}$) for the top (resp. bottom) of the resonance region. For the moment we restrict our attention to the top of the resonance region. 
On $\tr=0$, the linearisation of the vector field in $\delta \vect x$-coordinates \eqref{eq:linearization_equilibria} has the form,
\[\dot {\delta \vect x} = 2\pi \begin{pmatrix}
\varepsilon S_0 & C\Gamma - \varepsilon SS_0\\
-\varepsilon C_0& - \varepsilon S_0
\end{pmatrix}\delta \vect x, \]
where $S_0=\sin 2\pi x_0$, $C_0=\cos 2\pi x_0$ and $\Gamma = \sqrt{1-\varepsilon^2S_0^2}$ (the positive root corresponds to the centres at the top of the resonance region). On centres the determinant is positive so 
\[ 4\pi^2\varepsilon \left (C C_0\Gamma -\varepsilon S_0( S  C_0 +  S_0)\right ) >0.\]
In particular, for $\varepsilon$ small enough we have $C_0>0$. Moreover, as the determinant is invariant under change of bases, we have
\[\omega^2 = 4\pi^2\varepsilon \left (C C_0\Gamma -\varepsilon S_0( S  C_0 +  S_0)\right ), \]
and to simplify notation we write 
\[\omega = 2\pi f_0, \hspace{1cm} f_0= \sqrt{\varepsilon \left (C C_0\Gamma -\varepsilon S_0( S  C_0 +  S_0)\right )}. \]
The equations for the coordinate change \eqref{eq:Lyapunov_coef:coord_change} are given by,
\[2\pi \begin{pmatrix}
\varepsilon S_0 a - \varepsilon C_0 b & (C\Gamma -\varepsilon S S_0) a - \varepsilon S_0 b\\
\varepsilon S_0 c -\varepsilon C_0 d &
(C\Gamma -\varepsilon S S_0) c - \varepsilon S_0 d 
\end{pmatrix}
=
2\pi \begin{pmatrix}
-f_0 c & -f_0 d \\
f_0 a & f_0 b 
\end{pmatrix}. \]
As the equations are linearly dependent we have two degrees of freedom (scaling and rotation) in the choice of solution. We have chosen
\[a=f_0, \hspace{22pt} b=f_0, \hspace{22pt}c=\varepsilon(C_0-S_0),\hspace{22pt}d=-C\Gamma + \varepsilon S_0(S+1).\]
The determinant of this change of variable is given by,
\[\Delta =  -f_0(C\Gamma +\varepsilon C_0-\varepsilon S_0(S+2) ).  \]
Taylor expanding the vector field \eqref{eq:principal} about the centre $\vect x_0$ to third order in $\delta \vect x$ yields,
\begin{align*}
    \dot {\delta x} &= 2\pi \big ( \varepsilon S_0 \delta x+(C\Gamma - \varepsilon SS_0) \delta y \big ) +2\pi ^2\left ( \varepsilon C_0 \delta x ^2 + (S\Gamma +\varepsilon C S_0)\delta y^2 \right ) \\
    &\hspace{13pt}-\frac{4}{3}\pi ^3\left ( \varepsilon S_0 \delta x^3+(C\Gamma - \varepsilon SS_0) \delta y ^3 \right )\\
    \dot {\delta y} &=-2\pi \big ( \varepsilon C_0 \delta x + \varepsilon S_0 \delta y \big ) + 2\pi ^2 \left( \varepsilon S_0 \delta x^2 + \Gamma \delta y^2 \right)+ \frac{4}{3}\pi ^3 \left ( \varepsilon C_0 \delta x^3 + \varepsilon S_0 \delta y^3 \right ).
\end{align*}

We have already designed the change of variables to put the linear part in the form of \eqref{eq:Lyapunov_coef:linear_normal_form}, so we write the transformed vector field to third order as
\[\dot{\tilde {\vect x}}= \begin{pmatrix}
0& -\omega\\ \omega & 0
\end{pmatrix}\tilde{\vect x}
+
\begin{pmatrix}
P\\Q
\end{pmatrix}.
\]
We compute, 
\begin{align*}
    P = &\frac{2\pi ^2}{\Delta^2}\Big( (a\varepsilon C_0+b\varepsilon S_0)(d\tilde x- b \tilde y)^2+ (a (S\Gamma + \varepsilon C S_0)+ b\Gamma)(-c\tilde x + a \tilde y)^2\Big )\\
    &+\frac{4\pi ^3}{3\Delta^3}\Big ( (-a\varepsilon S_0 +b\varepsilon C_0)(d\tilde x-b\tilde y)^3+(-a(C\Gamma -\varepsilon S S_0)+b\varepsilon S_0)(-c\tilde x + a \tilde y)^3 \Big) 
\end{align*}
and,
\begin{align*}
    Q = &\frac{2\pi ^2}{\Delta^2}\Big( (c\varepsilon C_0+ d\varepsilon S_0)(d\tilde x- b\tilde y)^2 +(c(S\Gamma +\varepsilon C S_0)+d\Gamma)(-c\tilde x+a \tilde y )^2\Big )\\
    &+\frac{4\pi ^3}{3\Delta^3}\Big ( (-c\varepsilon S_0+d\varepsilon C_0)(d\tilde x- b\tilde y)^3 + (-c(C\Gamma - \varepsilon S S_0)+d\varepsilon S_0)(-c\tilde x + a \tilde y)^3  \Big) 
\end{align*}
%NOTE: we don't put epsilon in front of a,b,c,d as it helps to understand how we computed it. 

Following \cite{Baesens_2018}, which uses the notation in \cite{Nonlinear_oscillations_1983} with a different multiplicative factor, we compute the first Lyapunov coefficient as
\begin{align*}
    l_1 = &\frac{1}{8\omega}\big (P_{xxx}+P_{xyy}+Q_{xxy}+Q_{yyy}\big)\\
    &+\frac{1}{8\omega^2}\big (P_{xy}(P_{xx}+P_{yy})-Q_{xy}(Q_{xx}+Q_{yy})-P_{xx}Q_{xx}+P_{yy}Q_{yy}\big),
\end{align*}
where we have dropped the tildes from $x$ and $y$ and the derivatives are evaluated at 0. We will only need to compute this coefficient to first order in $\varepsilon$, so we use the following approximations, 
\begin{alignat*}{3}
\Gamma &\sim 1, \hspace{30pt}\ &f_0 &= a =b\sim \sqrt{\varepsilon C C_0},\\
\Delta &\sim -\sqrt{\varepsilon C^3 C_0},\hspace{30pt} &c&=\varepsilon(C_0-S_0),\hspace{30pt} d\sim-C.
\end{alignat*}
The $P$ terms are
\begin{alignat*}{3}
P_{xx} &\sim \frac{4\pi ^2 \sqrt{\varepsilon}}{\sqrt{CC_0} }\big ( C_0+S_0 \big)\hspace{30pt} 
&P_{xy} &\sim\frac{4\pi ^2 \varepsilon}{C^2}\big ( (C-S-1)C_0+ (C+S+1)S_0 \big)\\
P_{yy}&\sim\frac{4\pi ^2 \sqrt{\varepsilon C_0}}{\sqrt{C^3}}\big(S+1 \big) \hspace{30pt} &
P_{xxx}&\sim\frac{8\pi^3}{CC_0}\big (C_0-S_0\big) \hspace{30pt} 
P_{xyy}=0.
\end{alignat*}
% NOTE:  some of this coefficents could happen to be 0 (eg C0-S0) and then the higher order terms would play a role, so the sim notation is not fully justified. In this case however, they would be higher order and wouldn't play a role, so the computation of the lyapunov coef is correct.
The $Q$ terms are
\begin{alignat*}{3}
Q_{xx}&\sim-\frac{4\pi^2 S_0}{C_0} \hspace{30pt} 
& Q_{xy}&\sim\frac{4\pi^2\sqrt{\varepsilon}}{\sqrt{C^3 C_0}}\big (C_0-(C+1)S_0 \big )\\
 Q_{yy}&\sim-\frac{4\pi^2}{C} \hspace{30pt} &
Q_{yyy}& = o(\varepsilon) \hspace{30pt} 
Q_{xxy}\sim-\frac{8\pi^3}{C}.
\end{alignat*}
Finally we compute, 
\begin{equation}
\label{eq:Lyapunov_coef:l1_first_order}
    l_1\sim -\frac{\pi^2}{2\sqrt{\varepsilon C^7 C_0^5}}  \big( C S_0^2+S_0C_0+SC_0^2  \big).
\end{equation}
 So we have that $l_1<0$ as long as the bilinear form (on $S_0$ and $C_0$) inside the parentheses is positive definite. The condition of positive definiteness reduces to, 
\[CS-\frac{1}{4}>0,\]
or equivalently  $\phi \in (\frac{1}{24},\frac{5}{24})$ (assuming $\phi\in (-\frac{1}{4},\frac{1}{4}) $ which is needed to avoid degenerate cases of the ellipse $\psi$).

Recall that we have restricted our attention to the top of the resonance region. On the bottom, essentially the same computation yields \eqref{eq:Lyapunov_coef:l1_first_order} with opposite sign and absolute value inside the square root as in this case $C_0<0$. Thus, the first Lyapunov coefficient is non-zero and the bifurcation is generic.  
% NOTE notation a bit confusing because S_0 corresponds to s used in the previous section. However here we already use abcd for change variable. Moreover, we need x to make change of variables centererd at center which needs another name. 

\section{Average of \texorpdfstring{$f$}{f} inside chc is bigger than \texorpdfstring{$f(x_{\sad})$}{f(xsad)}}
\label{append:avarage_f}

For $\rho= \pm 1$  the chc degenerates to a $B$ point so there is nothing to prove. 
% EXPAND: Note that if we consider B point as chc then  average same as value at saddle (need stric inequality). This is not a real problem as B point is not consider as chc nor ns so there is no J point which is what we wanted to show. 
%\expand{For $rho$ bigger no equilibria so also trivial no $J$ point.}
We proceed to show that for $\rho\in [\rho_0,1)$ the average of $f$ defined in (\ref{eq:f}) over the region bounded by the chc (weighted by $e^{-2\pi x/C}$) is strictly bigger than $f(x_{\sad})$. A similar approach can be taken for $\rho \in (-1, \rho_0)$.

 We first show that the minimum of $g$ in the chc is attained at the saddle. Recall from  section~\ref{sec:reversible_approx}  that  $E_{\sad}\geq 0$ as $\rho \geq \rho_0$. Moreover, if $I_\rho$ is the interval of $x$ where the chc is defined, $x_{\sad}$ is its right endpoint and $x_{\cen}\in I_\rho$. Thus, $e^{-2\pi x /C }\geq e^{-2\pi x_{\sad} /C }$ for $x\in I_\rho$, so that the expression defining the chc,  $H(x, \eta) = E_{\sad}$, implies that
\[\eta^2 - \frac{1}{4\pi}g(x)\leq - \frac{1}{4\pi}g(x_{\sad}).\]
%\expand{$E_{\sad} = H(x_{\sad},0)>0$. Now evaluate and the equality thay you get use that $e^{-2\pi x/C}\geq e^{-2\pi x_{\sad}/C}$ together with the fact that the whole quantity is positive.   }
In particular, we have $g(x) \geq g(x_{\sad})$ as desired\footnote{For $\phi =0$, as was the case in \cite{Baesens_2018}, we are done as $f=g$.  }. 
%\expand{When $\rho <\rho_0$, $E_{\sad}<0$ but the $x_{\sad}$ is the left endpoint and we get the same result. }

Going back to $f$ it will be convenient to write
\[f(x)=\sqrt{(C-S)^2+1}\,\cos 2\pi (x-x_f),\]
where $x_f=\frac{\arctan (C-S)}{2\pi}$, again with the standard arctangent determination. Note that $x_f\leq x_g$ and they represent the phase shift of $f$ and $g$ respectively, so the minimum of $f$ in $I_\rho$ is also at $x_{\sad}$, as long as $I_\rho$ does not contain a local minimum of $f$.
% EXPAND: f is shifted to the left of g (as graphs) and not much so minimu is the right  of interval if not contains local minimum (the one you will reach, no chance of reaching other side as we would have had minimum in g also then). Note start (B point) interval is a point so no contains minimum. 
Finding when $x_{\sad}$ meets the first local minimum $x_f+\frac{1}{2}$,
%EXPAND: Need to meet this minimum first because we have x_cen (which is in interval of def) between x_f-1/2 and x_f+1/2. Moreover, if we contain minimum x_f-1/2 then we also contain minimum x_g-1/2 which is a contradiction as minimum of g is to the right hand side (and at most we are defined in length 1. 
one deduces that for $\rho\in [\rho_s,1)$, where
\[\rho_{s}   =-\frac{C-S}{\sqrt{2-2CS}},\]
%EXPAND: arcsin(rho)/2 pi = x_sad = x_f+1/2 = arctan(C-S)/2pi +1/2, multiplying by 2pi and then applying sin both sides we get formula for rho (use wolfram)
we have $f(x)\geq f(x_{\sad})$ in the chc, and thus the average is strictly bigger.
%\expand{Note that need strict inequality so that the two curves not meet. To deduce strict here note that f is never constant and interval is never a point as we have discarted B points.}

Consider now $\rho \in [\rho_0,\rho_s)$. We need to prove that,
\begin{equation*}
    \int_{\bar\gamma} (f(x)-f(x_{\sad}))e^{-2\pi x/C} \,\dd x \,\dd \eta = 2\int_{I_\rho} (f(x)-f(x_{\sad}))e^{-2\pi x/C} \eta(x)\,\dd x ,
\end{equation*}
is positive, where $\eta(x) $ is in the part of the chc   $\gamma$ above the $x$-axis. We do this by bounding the positive and negative contributions of the integrand separately. As $f$ is symmetric about the minimum $x_f+\frac{1}{2}$, we have that the integrand is negative exactly in $(x_{\midd},x_{\sad})$ where $x_{\midd}= 2x_f + 1-x_{\sad}$. Denote by $\eta^-_{\max}$ the supremum of $\eta$ in this interval. 
%NOTE: cannot have points outside the interval being negative as then we would get contradiction with x_sad being minimum of g (also it would need to be defined for atleast a unit of x which is impossible)
For the positive contribution it is convenient to only consider the interval $[x_{\cen},x_{\midd}]$ and denote by $\eta^+_{\min}$ the minimum $\eta  $ in it. 
%Note: we do this to take eta out and thus be able to integrate
Note that $\eta^-_{\max}\geq \eta^+_{\min}$.
%EXPAND: Needed latter this fact
So we find that it is enough to show that
\[\eta^+_{\min}\int_{x_{\cen}}^{x_{\midd}} (f(x)-f(x_{\sad}))e^{-2\pi x/C} \,\dd x>\eta^-_{\max}\int_{x_{\midd}}^{x_{\sad}} (f(x_{\sad})-f(x))e^{-2\pi x/C} \,\dd x,\]
which by integrating is equivalent to,
\[
    \eta^+_{\min}F(x_{\cen})>\eta^-_{\max}F(x_{\sad})-(\eta^-_{\max}-\eta^+_{\min})F(x_{\midd}),
\]
where 
% Note: integral of right hand side in equation to times above
\[F(x)=\frac{C e^{-2\pi x/C}}{2 \pi (1+C^2)}\left ( f(x)+\frac{C}{2 \pi} f'(x)-(1+C^2)f(x_{\sad})\right).\]
From the symmetry of $f$ about the minimum $x_f+\frac{1}{2}$ we have that $f(x_{\midd})=f(x_{\sad})$ and $f'(x_{\midd})=-f'(x_{\sad})$. Using this together with the fact that $\sin 2 \pi x_{\sad}=\sin 2 \pi x_{\cen}=\rho$ one can show that the term inside the parenthesis of $F(x_{\cen})$ is positive 
% EXPAND: it is needed as otherwise should do opposite bound
and bigger than the ones inside $F(x_{\sad})$ and $-F(x_{\midd})$.
% EXPAND: Should understand that minus goes inside parenthesis in $-F(x_{\midd})$. Also no need this terms to be positive. 
% % EXPAND:
% If we denote by $\tilde F$ the term inside the parenthesis of $F$ then we have
% \[\tilde F(x_{\cen}) = -C(1+C^2-CS)\rho +(2+2C^2-CS) \sqrt{1-\rho^2},\]
% \[\tilde F(x_{\sad}) = -C(1+C^2-CS)\rho +CS \sqrt{1-\rho^2},\]
% \[-\tilde F(x_{\mid}) = -C(1-C^2+CS)\rho -C(2C-S) \sqrt{1-\rho^2}.\]
% Then we can see that the coefficients in front of $\sqrt{1-\rho^2}$ of $\tilde F(x_{\cen})$, $\tilde F(x_{\cen})-\tilde F(x_{\sad})$ and $\tilde F(x_{\cen})+\tilde F(x_{\mid})$, are positive. So all three expressions are concave in $\rho$, and we can find their minimum evaluating $\rho$ at end points of the interval $[\rho_0, \rho_s]$. It is not hard to check that all this values are positive. 
Thus, it is enough to show that,
\begin{equation}
\eta^+_{\min}e^{-2\pi x_{\cen}/C}>\eta^-_{\max}e^{-2\pi x_{\sad}/C} +(\eta^-_{\max}-\eta^+_{\min})e^{-2\pi x_{\midd}/C}.    
\label{eq:eta^+_min}
\end{equation}
%\expand{Let $\tilde F(x_{\cen})$ be term inside parenthesis. As it is positive, \[\eta^+_{\min}e^{-2\pi x_{\cen}/C}>\eta^-_{\max}e^{-2\pi x_{\sad}/C} +(\eta^-_{\max}-\eta^+_{\min})e^{-2\pi x_{\midd}/C}\]   is equivalent to    \[\eta^+_{\min}e^{-2\pi x_{\cen}/C}\tilde F(x_{\cen})>\eta^-_{\max}e^{-2\pi x_{\sad}/C}\tilde F(x_{\cen}) +(\eta^-_{\max}-\eta^+_{\min})e^{-2\pi x_{\midd}/C}\tilde F(x_{\cen}).  \]     As $\tilde F(x_{\cen})>\tilde F(x_{\sad}),-\tilde F(x_{\midd})$ and  $\eta^-_{\max},(\eta^-_{\max}-\eta^+_{\min})>0$, the right hand side above is an upper bound to     \[\eta^-_{\max}F(x_{\sad})-(\eta^-_{\max}-\eta^+_{\min})F(x_{\midd}).\] Putting the last two inline equations together we get the desired inequality. }
Now looking at the nullclines of the Hamiltonian system \eqref{eq:reversible_approx}, we find 
that $\eta (x)$ has a unique relative extremum $\tilde x$, which is a maximum with value $\sqrt{\sin 2\pi \tilde x - \rho}/(\pi \sqrt{2})$.
%EXPAND: in maximum we have n'=0 in the ODE so we get expression. Local extremum is unique because the n'=0 curve is contractible, so in n>0 can only have one intersection with our chc. 
Thus, $\eta_{\max}^-$ is bounded by $\pi^{-1}$, and $\eta_{\min}^+$ is attained either at $x_{\cen}$ or at $x_{\midd}$. In the latter case, we have $\eta_{\max}^-=\eta_{\min}^+$,
%EXPAND: Because in this case we necessarly have maximum eta in between x_min and x_mid. So x_mid has to be maximum for interval [x_mid, x_sad] (where eta is decreasing). 
so \eqref{eq:eta^+_min} becomes, $e^{-2\pi (x_{\cen}-x_{\sad})/C}>1$ which is satisfied as $x_{\cen}<x_{\sad}$.

If $\eta_{\min}^+$ is attained at $x_{\cen}$ we have, $H(x_{\cen},\eta_{\min}^+)=E_{\sad}$, and thus
%\expand{Also use the sinus at saddle and center are $\rho$ but already said this recently (from this we also get value of cosinus so value of g)}
\[(\eta_{\min}^+)^2=\frac{C}{2\pi^2(1+C^2)}\left (\sqrt{1-\rho^2}-C\rho+e^{-2\pi(x_{\sad}- x_{\cen})/C}(\sqrt{1-\rho^2}+C\rho)\right).\]
The minimum of the first and third term inside the parenthesis are attained at $\rho_0$, whereas for the second term it is attained at  $\rho_s$, so we get
%\expand{The second term is clearly decreasing so take right hand limit, $\rho_s$. The first term is concave so minimum at one of the boundaries, one can check (not completely trivial) that at $\rho_0$ value always smaller than at $\rho_s$. For third term both summands attain minimum at $\rho_0$.}
\[(\eta_{\min}^+)^2>\frac{-C^2}{2\pi^2(1+C^2)}\left ( \rho_0+\rho_s \right).\]
%\expand{\[\sqrt{1-\rho_0^2}= \sqrt{C^2/(1+C^2)}= -C\rho_0\] So we get first term and get 0 for third term. }
For $\rho>\rho_0$ we also have   $x_{\cen}=\frac{1}{2}-x_{\sad}$, so that $x_{\midd}-x_{\cen}=2x_f+\frac{1}{2}$. Then, if in \eqref{eq:eta^+_min} we bound the difference of $\eta$'s by $\eta_{\max}^-$ which in turn we bound by $\pi^{-1}$ and we also bound the exponentials of the right hand side by $e^{-2\pi x_{\midd}/C}$, we find that it is enough to show, 
 \[\frac{C}{\sqrt{2}\,\pi\sqrt{1+C^2}}\sqrt{- \rho_0-\rho_s }\, e^{2\pi (2x_f+\frac{1}{2})/C}>\frac{2}{\pi }.\]
We note that $C$, $C^2$, $\rho_0$, $\rho_s$  are monotone in $\phi$, so if we bound them by their value at $\frac{1}{24}$ or $\frac{5}{24}$ and move them to the other side we get, 
\[ e^{2\pi (2x_f+\frac{1}{2})/C}>\frac{24 + 4\sqrt{3}}{\sqrt{8 - 6 \sqrt{3} + 4 \sqrt{30 - 17 \sqrt{3}}}}.\]
 It is not hard to check that this inequality is true for $\phi \in (\frac{1}{24},\frac{5}{24})$.

\section{Integral of \texorpdfstring{$\tau$}{tau} in a neighbourhood  of the saddle}
\label{append:integral_trace_saddle}
Let $(x(s),\eta(s)) $ parameterise the part of the  orbit $\gamma$ in $U$, so that $(x(0),\eta(0))$ and $(x(T),\eta(T))$ are at the boundary of $U$. Denote by $h$ the smooth conjugacy that flattens the coordinates and by  $(u,v)$ the flattened coordinates so that there exist smooth $F_1, F_2$ such that, 
\begin{equation}
\begin{split}
    \primebullet u &= -u F_1(u,v)\\
    \primebullet v &= v F_2(u,v),
\end{split}
\label{eq:flattened_dynamics}
\end{equation}
with $F_1(0,0),F_2(0,0)>0$ and around\footnote{That is,  $F_i(0,0)=O(1)$ and $1=O(F_i(0,0))$ for $i=1,2$.} 1.
%\expand{Around 1 as egenvalues have been computed and are around 1. }
By shrinking $U$ independently of $\varepsilon$ we can assume that $F_1,F_2$ are around 1 in the whole of $h(U)$.
%\expand{Should check that this shrinking is independent of $\varepsilon$ which is is as eigenvalues constant to first order.}
Then, if $(u(s),v(s))$ is the flow of \eqref{eq:flattened_dynamics} starting at $h(x(0),\eta(0))$ we have $h(x(s),\eta(s))=(u(s),v(s))$. 

Here we treat $\tau$ from \eqref{eq:dot_A} as a function, so we write   $\tau (x,y) $, and we want to compute 
\[\int_0^T\tau (x(s), \eta(s))  \,\dd s = \int_0^T(\tau \circ h^{-1})(u(s), v(s))  \,\dd s \]
We define a smooth function $\tilde \tau (u,v)$ by, 
\begin{equation*}
    \sqrt{\varepsilon} \:\tilde \tau (u,v) =  (\tau \circ h^{-1}) (u,v),
\end{equation*}
so that we want to find $\sqrt{\varepsilon}\int_0^T\tilde \tau(u(s),v(s)) \, \dd s $. Note that from \eqref{eq:dot_A} it is clear that $\tilde \tau$ and its derivatives are $O(1)$.
%\expand{Note that by this we mean (what will use when differentiating last equation respect parameters inside integral) derivative respect parameters of $\tilde \tau_u$, so second order derivatives of $\tilde \tau$.}
Then, by shrinking $U$ further so that $h(U)$ is convex  we have,
%\expand{Need to restrict here as $\tilde \tau $ only defined where $h^{-1}$ defined. More formally we take $V\subset h(U)$ convex and transverse to flow (at least for $v$ small which corresponds to $\eta $ close to $\eta_s$). Then we redefine $U:=h^{-1}(V)$ which transverse as $h$ smooth, but maybe not convex (doesn't matter). }
\[\tilde \tau(u,v)=\tilde \tau(0,0) + u\tilde \tau_u (c(u,v))+v\tilde \tau_v (c(u,v)) \]
by the mean value theorem,  where $\tilde \tau_u \circ c$, $\tilde \tau_v \circ c$ are smooth. 
%\expand{Smooth as result comes from 1 variable result and there we have $\tilde \tau_u(c(u)) = \frac{\tilde \tau (u) -\tilde \tau (0) }{u} $ which is clearly smooth for $u\not =0$ and at 0 smooth as all limits behave well as $\tilde \tau$ smooth. Note however that need smoothnes respect parameters really. But same approach works just by introducing more variables in the function $\tilde \tau$ (and thinking segment considered doesn't change parameters when applying mean value theorem).   }
Thus, if we let $\tilde \gamma = h(\gamma\cap U)$ and use \eqref{eq:flattened_dynamics} to change variable of integration, we have
\[\sqrt{\varepsilon}\int_0^T\tilde \tau(u(s),v(s)) \, \dd s =\sqrt{\varepsilon}\tilde\tau (0, 0) T+\sqrt{\varepsilon}\left (-\int_{\tilde \gamma} \frac{\tilde\tau_u (c (u, v))}{F_1(u,v)}  \,\dd u+\int_{\tilde \gamma} \frac{\tilde\tau_v (c (u, v))}{F_2(u,v)}  \,\dd v\right ) .\]
%\expand{Here by v inside the frist integral we mean v(u) corresponding to values in the orbit.}
The integrals are bounded as $F_1, F_2$ are around 1 so the second summand is $O(\sqrt{\varepsilon})$.  Moreover,  the order is preserved upon differentiation with respect to parameters by differentiating inside the integral sign. Finally, as conjugacies send equilibria to equilibria, $\sqrt{\varepsilon}\tilde{\tau}(0,0) T = T \tau_{\sad}$. So we have shown \eqref{eq:int_tau} as desired.

\printbibliography
\end{document}